\documentclass[sn-mathphys,Numbered]{sn-jnl}% Math and Physical Sciences Reference Style
\usepackage{graphicx}%
\usepackage{multirow}%
\usepackage{amsmath,amssymb,amsfonts}%
\usepackage{amsthm}%
\usepackage{mathrsfs}%
\usepackage[title]{appendix}%
\usepackage{xcolor}%
\usepackage{textcomp}%
\usepackage{manyfoot}%
\usepackage{booktabs}%
\usepackage{algorithm2e}
\usepackage{algpseudocode}%
\usepackage{listings}%
\usepackage{caption}
\usepackage{subcaption}
\usepackage{hyperref}
\renewcommand{\leq}{\leqslant}
\renewcommand{\geq}{\geqslant}

\renewcommand{\ell}{j}

\renewcommand{\circ}{*}

\lstdefinelanguage{Julia}{
	keywords={function, end, if, else, elseif, while, for, return, break, continue, begin, type, struct, mutable, using, import, export, const, let, global, local, true, false, abstract, primitive, module, baremodule, where, do},
	sensitive=true,
	comment=[l]{\#},
	morestring=[b]",
	morestring=[b]',
}

\begin{document}

%\title[SIMD-vectorized implementation of high order implicit Runge-Kutta integrators for non-stiff ODEs]{SIMD-vectorized implementation of high order implicit Runge-Kutta integrators for non-stiff ODEs}
\title[SIMD-vectorized implicit symplectic integrators can outperform explicit ones]{SIMD-vectorized implicit symplectic integrators can outperform explicit symplectic ones}


\author[1]{\fnm{Mikel} \sur{Antoñana} (Corresponding author)}\email{mikel.antonana@ehu.eus}

\author[1]{\fnm{Joseba} \sur{ Makazaga}}\email{joseba.makazaga@ehu.eus}
\equalcont{These authors contributed equally to this work.}

\author[1]{\fnm{Ander} \sur{Murua}}\email{ander.murua@ehu.eus}
\equalcont{These authors contributed equally to this work.}

\affil[1]{\orgdiv{Computer Science and Artificial Intelligence Department}, \orgname{UPV/EHU (University of the Basque Country)}, \orgaddress{\street{} \city{Donostia}\postcode{}, \state{} \country{Spain}}}


\abstract{
Implicit Runge-Kutta schemes based on collocation with Gauss-Legendre nodes (IRKGL) possess attractive theoretical properties for long-term integration of Hamiltonian systems—they are symmetric, symplectic, and super-convergent—yet they have traditionally been considered less practical than explicit symplectic integrators. In this work, we challenge this conventional wisdom by introducing a stage-wise SIMD-vectorization approach that enables IRKGL schemes to outperform state-of-the-art explicit symplectic integrators for high-precision computations in double-precision floating-point arithmetic. Our approach reformulates the fixed-point iteration in terms of $s$-vectors (where $s$ is the number of stages) to explicitly exploit Single Instruction Multiple Data (SIMD) capabilities of modern processors, significantly reducing computational overhead and enabling parallel evaluation of the right-hand side function across all stages. Additionally, we present a reformulation of IRKGL schemes for second-order ODEs that ensures exact symplecticity at the level of double-precision floating-point arithmetic, extending our previous approach for first-order systems to this setting. We demonstrate these ideas through IRKGL16, a Julia implementation of the 16th-order 8-stage IRKGL scheme that performs vectorization seamlessly and transparently to the user. Numerical experiments on several Hamiltonian problems with separable structure confirm the effectiveness of our approach.}
	



\maketitle
\pagebreak

\section{Introduction}

The family of $s$-stage implicit Runge-Kutta schemes based on collocation with Gauss-Legendre nodes (IRKGL) are known to be symmetric, symplectic and super-convergent of order $2s$, making them theoretically well-suited for high-precision long-term numerical integration of Hamiltonian systems~\cite{Sanz1994,Hairer2006}. However, despite these attractive theoretical properties, IRKGL schemes have traditionally been considered less practical than explicit symplectic integrators for Hamiltonian problems that admit a separable structure. The main reasons for this preference include: (i) the need for careful implementation to avoid linear drift of energy error~\cite{Hairer2008,antonana2017}, (ii) the higher computational overhead inherent to implicit methods, and (iii) the typically larger number of right-hand side evaluations required to achieve comparable accuracy.

For Hamiltonian systems with suitable structure---where the Hamiltonian function can be decomposed into the sum of two or more exactly solvable components---optimized high-order explicit symplectic integrators derived through splitting and composition techniques~\cite{Blanes2024,Blanes2016,Hairer2006,Leimkuhler2004} have thus become the methods of choice for high-precision integration. These methods avoid the complications associated with implicit schemes while maintaining the crucial symplectic property.

In this work, we challenge this conventional wisdom by showing that IRKGL schemes can outperform state-of-the-art explicit symplectic integrators for high-precision computations in double-precision floating-point arithmetic. Our main contribution is a novel approach to formulating and implementing the fixed-point iterations of IRKGL schemes: we reformulate the iterations in terms of $s$-dimensional vectors (where $s$ is the number of stages) and design the algorithm to explicitly exploit Single Instruction Multiple Data (SIMD) capabilities of modern processors. 

This stage-wise SIMD-vectorization approach represents a departure from conventional strategies for reducing overhead in implicit integrators. Traditionally, SIMD-vectorization is achieved by structuring inner loops in a way that allows the compiler to automatically vectorize them. In contrast, our approach explicitly builds SIMD operations into the algorithm design itself. This has two crucial advantages: it significantly reduces the computational overhead of the implicit method, and it enables parallel evaluation of the right-hand side function across all $s$ stages simultaneously (provided the function satisfies certain requirements for vectorization). 
%It is well established~\cite{Hairer2006} that implicit Runge-Kutta methods can benefit substantially from parallelism; our stage-wise vectorization exploits SIMD capabilities where vector operations can be performed at costs comparable to scalar operations.

To ensure exact symplecticity at the level of double-precision floating-point arithmetic, we address the linear drift in energy error that can affect symplectic IRK schemes~\cite{Hairer2008}. For first-order ODEs, we adopt the approach presented in~\cite{antonana2017}. Additionally, we generalize this approach to the implicit Runge-Kutta-Nyström (IRKN) formulation for second-order ODEs, which is the natural framework for many Hamiltonian systems arising from classical mechanics.

We demonstrate these ideas through a concrete implementation, IRKGL16, which applies our stage-wise SIMD-vectorization strategy to the $16$th-order $8$-stage IRKGL scheme (i.e., $s=8$, order $2s=16$). The implementation, written in Julia using the \verb|SIMD| package~\cite{SIMD}, performs vectorization automatically and seamlessly, remaining transparent to the user. While IRKGL16 includes both constant step-size and adaptive step-size modes, in this work we focus on the constant step-size mode, as recommended for long-term symplectic integration of Hamiltonian systems~\cite{Sanz1994,Hairer2006,Leimkuhler2004}. We have also implemented IRKGL methods with $s=2$, $s=4$, and $s=6$ stages, which are available in the IRKGL16 code, though our numerical experiments (see Subsection~\ref{ss:IRKGL_comp}) demonstrate that the $8$-stage method achieves the best performance for high-precision computations.

The code is designed for non-stiff systems of ODEs (not necessarily Hamiltonian). While it can handle arbitrary precision arithmetic, the vectorized implementation is available for single- and double-precision floating-point arithmetic. Maximum efficiency through full vectorization is achieved when the right-hand side function is implemented as a generic Julia function using arithmetic operations and elementary functions supported by the SIMD package.

Through numerical experiments on several representative Hamiltonian problems with separable structure, we show that our stage-wise SIMD-vectorization approach enables IRKGL16 to outperform state-of-the-art explicit symplectic integrators for high-precision integration tasks in double-precision arithmetic.

The paper is organized as follows: Section 2 describes the reformulation of IRKGL schemes for both first-order and second-order ODEs to ensure exact symplecticity and our stage-wise SIMD-vectorization strategy; Section 3 presents implementation details of IRKGL16; Section 4 reviews state-of-the-art high-order explicit symplectic integrators; Section 5 reports numerical experiments; and Section 6 presents the conclusions.

\section{Vectorized fixed-point implementation of symplectic IRK schemes}

\subsection{Implicit RK methods}

In this work, we focus on the numerical integration of autonomous Hamiltonian systems of the form
\begin{equation}
\label{eq:Ham}
\frac{d y}{dt} = J^{-1} \nabla H(y), \quad y \in \mathbb{R}^{2d},
\end{equation}
where $H: \mathbb{R}^{2d} \to \mathbb{R}$ is the Hamiltonian function and $J$ is the standard symplectic matrix
\[
J = \begin{pmatrix}
0 & I_d \\
-I_d & 0
\end{pmatrix}.
\]
Runge-Kutta (RK) methods are widely used one-step integrators for numerically solving ODEs. These methods are specified by the coefficients in the Butcher tableau:
\begin{equation*}
\begin{array}{c|ccc}
c_1 & a_{11} & \dots & a_{1s} \\
\vdots & \vdots & \ddots & \vdots \\
c_s & a_{s1} & \dots & a_{ss} \\
\hline
& b_{1} & \dots & b_{s} \\
\end{array}
\end{equation*}
Although our focus is on Hamiltonian systems, RK methods are applicable to more general ODE problems of the more general form:
\begin{equation}
\label{eq:ivp}
\frac{d}{dt}y=f(t,y),\quad  y(t_0)=y_0,
\end{equation}
where $f: \mathbb{R}^{D+1} \to \mathbb{R}^D$ is a sufficiently smooth function and $y_0 \in \mathbb{R}^D$.

The RK approximations $y_n \approx y(t_n)$ to the solution $y(t)$ of \eqref{eq:ivp} at times $t=t_n = t_0 + nh$ for $n = 1, 2, 3, \dots$ are given by:
\begin{equation}
\label{eq:yn}
y_{n} = y_{n-1} + h \sum_{i=1}^s b_i f(t_{n-1} + c_i h, Y_{n,i}),
\end{equation}
where
\begin{equation}
\label{eq:Yi}
Y_{n,i} = y_{n-1} + h \sum_{j=1}^s a_{ij} f(t_{n-1} + c_j h, Y_{n,j}), \quad i = 1, \ldots, s.
\end{equation}
Here, each $Y_{n,i}$ is an approximation of the state vector $y(t)$ at the intermediate time $t=t_{n-1}+h c_i$ computed within each time-step. The vectors $b=(b_1,\ldots,b_s)$, $c=(c_1,\ldots,c_s)$, and the matrix $A=(a_{ij})$ from the Butcher tableau determine the specific RK method and its properties.

If the matrix $A$ is lower triangular %(potentially after reordering the indices), 
the Runge-Kutta method is explicit. Otherwise, it is an implicit Runge-Kutta (IRK) scheme, and the stage vectors $Y_{n,i}$ at each step are defined implicitly by the equation (\ref{eq:Yi}).

The primary challenge in implementing implicit Runge-Kutta (IRK) methods is the efficient solution of the nonlinear system for the internal stages. For non-stiff problems, fixed-point iteration is generally recommended~\cite{Sanz1994,Hairer2006}. 

\subsection{IRK methods of collocation type}

Among IRK methods, those based on collocation are particularly attractive.
In a collocation method based on the pairwise distinct nodes $c_1,\ldots,c_s \in [0,1]$,  
the approximations $y_n \approx y(t_n)$ to the solution $y(t)$ at times $t_n = t_0 + nh$ are obtained as follows: for each $n = 1, 2, 3, \dots$, consider the polynomial function $P_n:\mathbb{R} \to \mathbb{R}^D$ of degree $s$ uniquely determined by the following conditions,
\begin{itemize}
\item $P_n(t_{n-1}) = y_{n-1}$,
\item $P_n'(t_{n-1}+c_i h) = f(t_{n-1}+c_i h, P_n(t_{n-1}+c_i h))$ for $i=1,\ldots,s$.
\end{itemize}
then, set $y_n = P_n(t_n)$.

This collocation procedure is equivalent~\cite{Hairer2006} to the application of the IRK method (\ref{eq:yn})--(\ref{eq:Yi}) with coefficients $b_i,a_{ij}$ uniquely determined by the following conditions: 
\begin{equation}
\begin{split}
\sum_{j=1}^s a_{i j} c_j^{k-1} &=\frac{c_i^k}{k}, \quad i = 1,\ldots,s, \quad k=1,2,\ldots,s, \\
\sum_{j=1}^s b_j c_j^{k-1} &=\frac{1}{k}, \quad k=1,2,\ldots,s.
\end{split}
\end{equation}
In the specific case where $c_1,\ldots,c_s \in (0,1)$ are the shifted Gauss-Legendre nodes (defined by $c_i =(1+x_i)/2$, where $x_1,\ldots,x_s$ are the zeros of the Legendre polynomial of degree $s$) the resulting IRK scheme of collocation type, which we refer to as the $s$-stage IRKGL scheme,  achieves an order of accuracy $2s$~\cite{Hairer2006}.

In collocation-type IRK methods, it is standard practice to initialize the fixed point iteration with an estimate of the internal stages extrapolated from the previous time-step, that is, 
\begin{equation}
\label{eq:Yni0}
Y_{n,i} \approx P_{n-1}(t_{n-1} + h c_i), \quad i = 1,\ldots,s.
\end{equation}
It is not difficult to check that 
\begin{equation}
\label{eq:Pol0}
P_{n-1}(t_{n-1} + h c_i) = y_{n-1} + \sum_{j=1}^{s} \nu_{i j} L_{n-1,j}, \quad i = 1,\ldots,s,
\end{equation}
where the coefficients $\nu_{ij}$ are uniquely determined from the following equations
\begin{equation*}
\sum_{j=1}^s \nu_{ij} (c_i-1)^{k-1} = \frac{c_i^k}{k}, \quad i=1,\ldots,s, \quad k=1,\ldots,s.
\end{equation*}

\subsection{Application of IRK methods to second order differential equations}

Autonomous Hamiltonian systems with a Hamiltonian function of the form
\begin{equation}
\label{eq:sHam}
H(q,p) = \frac12\, p^T M^{-1} p + U(q)
\end{equation}
can be expressed as system of second order ODEs of the form
\begin{equation}
\label{eq:Hamode2}
\frac{d^2}{dt^2} q = M^{-1} \nabla U(q).
\end{equation}
Any RK scheme can be applied to systems of second-order differential equations $\frac{d^2}{dt^2} q = g(t,q)$ by rewriting it as a system of first-order differential equations,
\begin{equation}
\label{eq:ivp2}
\frac{dq}{dt}=v, \quad \frac{dv}{dt}=g(q,t),
\end{equation}
Obviously, the system (\ref{eq:ivp2}) is a specific instance of the general form
(\ref{eq:ivp}), with $D=2d$, $y=(y^1,\ldots,y^{2d})$, $q=(y^1,\ldots,y^{d})$, $v =(y^{d+1},\ldots,y^{2d})$, and $f(t,y) = (v,g(t,q))$.

It is well known that the implicit equations of the IRK scheme can be simplified in the case of second order ODE systems of the form \eqref{eq:ivp2} by rewriting it as an implicit Runge-Kutta-Nyst\o m (IRKN) method. More precisely, the approximations $y_n = (q_n,v_n)\approx (q(t_n),v(t_n))$ of the solution $y(t)=(q(t),v(t))$ of the  \eqref{eq:ivp2} supplemented with the initial conditions 
\begin{equation}
\label{eq:ic2}
 q(t_0)=q_0, \quad v(t_0)=v_0
\end{equation}
can be computed as follows: for $n=1,2,\ldots$
\begin{equation}
\label{eq:qnvn}
\begin{split}
q_{n} &= q_{n-1} + h\, v_{n-1} + h^2 \sum_{i=1}^s \beta_i\, g(t_{n-1} + c_i h, Q_{n,i}), \\
v_{n} &= v_{n-1} + h \sum_{i=1}^s b_i g(t_{n-1} + c_i h, Q_{n,i}),
\end{split}
\end{equation}
where the vectors $Q_{n,i}$, $i=1,\ldots,s$,  are implicitly defined by
\begin{equation}
\label{eq:Qi}
Q_{n,i} = q_{n-1} + h\, c_i \, v_{n-1} + h^2\, \sum_{j=1}^s \alpha_{ij}\, g(t_{n-1} + c_j h, Q_{n,j}), \quad i = 1, \ldots, s.
\end{equation}
Here,  $\beta_i = \sum_{k=1}^s b_k a_{k i}$, $c_i = \sum_{k=1}^s a_{i k}$,  and $\alpha_{i j} = \sum_{k=1}^s a_{i k} a_{k j}$.
This formulation gives rise to a more efficient fixed-point iteration, which, compared to the standard IRK formulation, approximately halves the required number of iterations for solving the implicit equations~\cite{Sanz1994,Hairer2006}.

\subsection{Symplectic IRK methods with floating-point coefficients}

\paragraph{Symplectic IRK methods}

The $t$-flow of the Hamiltonian system \eqref{eq:Ham}, denoted by $\varphi_t$, is symplectic, meaning that it preserves the symplectic form: for all $t$, the Jacobian matrix $D\varphi_t(y)$ satisfies
\[
\left(D\varphi_t(y)\right)^T J D\varphi_t(y) = J.
\]
It is therefore desirable that numerical integrators used to approximate such flows also preserve this structure, in order to ensure long-time stability and accurate qualitative behavior of the solution.

An implicit Runge-Kutta (IRK) method is symplectic if and only if~\cite{Sanz1994} its coefficients satisfy the condition
\begin{align}
\label{eq:sym}
b_i a_{ij} + b_j a_{ji} = b_i b_j, \quad 1 \leq i, j \leq s.
\end{align}
This condition ensures that the discrete flow defined by the IRK method is symplectic, and is satisfied in particular by the $s$-stage IRKGL schemes~\cite{Sanz1994,Leimkuhler2004}.

\paragraph{Symplectic IRKN methods}

An IRK method rewritten as an IRKN \eqref{eq:qnvn}--\eqref{eq:Qi} is symplectic if its coefficients satisfy the conditions~\cite{Suris1989,Sanz1994}
\begin{equation}
\label{eq:sym2}
\begin{split}
\beta_i &= b_i\, (1-c_i), \quad 1 \leq i \leq s, \\
b_i \,(\beta_j - \alpha_{ij}) &= b_j\,(\beta_i - \alpha_{ji}), \quad 1 \leq i, j \leq s.
\end{split}
\end{equation}

\paragraph{Reformulation for floating-point arithmetic}

As shown in~\cite{Hairer2008}, if the coefficients $b_i,a_{ij}$ of the IRK schemes (or $b_i, \beta_i, c_i, \alpha_{ij}$ of the IRKN formulation) are replaced by their machine number representatives, then the resulting implicit RK scheme typically fails to satisfy the symplecticity condition (\ref{eq:sym}). To ensure that the symplectic property is exactly preserved with double-precision floating-point coefficients, 
for first-order systems, we adopt the reformulation from~\cite{antonana2017}: 
\begin{align}
\label{eq:Yi2}
Y_{n,i}  &=y_{n-1}+ \sum^s_{j=1}{\mu_{ij}\,L_{n,j}, \quad  L_{n,i} = h b_i f(t_{n-1}+hc_i,Y_{n,i})}, \quad  i=1 ,\ldots, s, \\
\label{eq:yn2}
y_{n} &=y_{n-1}+\sum^s_{i=1} L_{n,i},
\end{align}   
where $b_i$ and $c_i$ are the machine-number representatives of the original coefficients of the IRK scheme, and 
\begin{align*}
\mu_{i j} &= \mathrm{fl}(a_{ij}/b_j), \quad 1 \leq j \leq i \leq s, \\
\mu_{i j} &= 1 -  \mu_{j i}, \quad 1 \leq i < j \leq s.
\end{align*}
Here,  $\mathrm{fl}(x)$ represents the machine-number representative of the real number $x$.
%The symplecticity conditions $\tilde \mu_{i j} +\tilde \mu_{j i} = 1$ hold for $s$-stage IRKGL schemes with $s\leq 8$~\cite{antonana2017}.
%This ensures that the symplecticity conditions $\tilde \mu_{i j} +\tilde \mu_{j i} = 1$ are exactly satisfied, because 
%$1/2 \leq \tilde \mu_{j i} \leq 2$ for $1 \leq j < i \leq s$.

For second-order systems, we generalize this approach as follows: We reformulate the IRKN form (\ref{eq:qnvn})--(\ref{eq:Qi}) of IRKGL schemes in terms of the coefficients $b_i$, $c_i$, and $\eta_{i j} = \alpha_{i j}/b_j$,
so that the second condition in (\ref{eq:sym2}) can be replaced by 
\begin{equation}
\label{eq:sympl_cond_2}
\eta_{i j} + c_j = \eta_{j i} + c_i.
\end{equation}

\smallskip

 Now, we define one step of the method as follows:
\begin{align}
\label{eq:Gi}
 G_{n,i} &= g(t_{n-1}+hc_i,Q_{n,i}), \quad  
 R_{n,i} = h \, b_i \, G_{n,i}, \quad  i=1 ,\ldots, s, \\
 \label{eq:Qi2}
Q_{n,i}  &=q_{n-1}+ h\, c_i \, v_{n-1} +  h \sum^s_{j=1}\eta_{ij}\,R_{n,j}, \quad  i=1 ,\ldots, s, \\
\label{eq:qnvn2}
v_{n} &=v_{n-1}+ \sum^s_{i=1} R_{n,i}, \quad
q_{n} =q_{n-1}+h\, v_{n} - h \sum^s_{i=1} c_i R_{n,i}
\end{align}   
where the coefficients $b_i,c_i$ are machine-number representatives of the original coefficients of the $s$-stage IRKGL scheme, and 
\begin{align*}
\eta_{i j} &= \mathrm{fl}(\alpha_{i j}/b_j), \quad 1 \leq j \leq i \leq s, \\
\eta_{i j} &= c_i - c_j - \eta_{j i}, \quad 1 \leq i < j \leq s.
\end{align*}
This guarantees that the symplecticity conditions \eqref{eq:sympl_cond_2} holds. 
We have verified numerically that the coefficients $\eta_{ij}$ computed by this procedure are exactly representable as machine numbers in double-precision arithmetic for $s \leq 16$.

\subsection{Stage-wise vectorization of fixed-point iteration}
\label{ss:vec_IRK}

In order to enable SIMD implementation, we reformulate the method using
$s$-vectors that group values across all stages. 

\paragraph*{Vectorized notation}

We introduce notation to facilitate the description of the proposed stage-wise vectorized implementation of the fixed-point iteration.

Let $\boldsymbol{b}$ and $\boldsymbol{c}$ denote the $s$-vectors $ (b_1, \ldots, b_s) $ and $ (c_1, \ldots, c_s) $, respectively. Similarly, for each $ i = 1, \ldots, s $, let $ \boldsymbol{\mu_i} $, $ \boldsymbol{\nu_i} $,  and $ \boldsymbol{\eta_i}$ represent the $ s $-vectors $ (\mu_{1i}, \ldots, \mu_{si}) $,   $ (\nu_{1i}, \ldots, \nu_{si}) $, and $ (\eta_{1i}, \ldots, \eta_{si}) $, respectively.

For an arbitrary $s$-vector $\boldsymbol{v}=(v_1,\ldots,v_s)$ and a real number $\lambda \in \mathbb{R}$, we denote by $\lambda + \boldsymbol{v}$ the $s$-vector $(\lambda +v_1,\ldots,\lambda + v_s)$. Given two arbitrary $s$-vectors $\boldsymbol{v}=(v_1,\ldots,v_s)$ and $\boldsymbol{w}=(w_1,\ldots,w_s)$, $\boldsymbol{v} \circ \boldsymbol{w}$ represents the componentwise product of two $s$-vectors defined by $(v_1 w_1, \ldots, v_s w_s)$. 
 We also use the notation $\mathrm{sum}(\boldsymbol{v}) = v_1+\cdots + v_s$, 
 %$\mathrm{dot}(\boldsymbol{v},\boldsymbol{w}) = v_1w_1+\cdots + v_sw_s$, 
 and $\|\boldsymbol{v}\|_{\infty} = \max(|v_1|,\ldots,|v_s|)$.

To distinguish $s$-vectors from collections of $D$ state variables, we denote the laters as 
\[
y=\left[
\begin{matrix}
y^1\\
\vdots\\
y^D
\end{matrix}
\right]
\]
 (or alternatively as $y=[y^1,\ldots,y^D]$) and refer to them as arrays of state variables.

For each $\ell \in \{1,\ldots,D\}$, we denote by $Y_{n,i}^{\ell}$, $L_{n,i}^{\ell}$, and $f^{\ell}(t_{n-1}+hc_i,Y_{n,i})$ the $\ell$th component of the arrays $Y_{n,i} \in \mathbb{R}^D$, $L_{n,i} \in \mathbb{R}^D$ and  $f(t_{n-1}+hc_i,Y_{n,i}) \in \mathbb{R}^D$ respectively.  
For each $\ell \in \{1,\ldots,D\}$, we consider the $s$-vectors $\boldsymbol{Y_{n}^{\ell}} = (Y_{n,1}^{\ell}, \ldots,Y_{n,s}^{\ell})$ and  $\boldsymbol{L_{n}^{\ell}} = (L_{n,1}^{\ell}, \ldots,L_{n,s}^{\ell})$.
 
Let $\boldsymbol{Y_n}$ and $\boldsymbol{F_n}$ represent the arrays of $s$-vectors $\boldsymbol{Y_n^\ell}$ and $\boldsymbol{F_n^\ell}$ 
\[
\boldsymbol{Y_n} = 
\left[
\begin{matrix}
\boldsymbol{Y_{n}^1}\\
\vdots\\
\boldsymbol{Y_{n}^D}
\end{matrix}
\right], \quad
\boldsymbol{F_n} =
\left[
\begin{matrix}
\boldsymbol{F_{n}^1}\\
\vdots\\
\boldsymbol{F_{n}^D}
\end{matrix}
\right].
\]
We define the \textbf{vectorized} function  $\boldsymbol{f}$ that given $t_{n-1} + h \boldsymbol{c}$ and $\boldsymbol{Y_n}$, returns the array $\boldsymbol{F_n}$  of $s$-vectors
\[
\boldsymbol{F_n^{\ell}}= (f^{\ell}(t_{n-1}+hc_1,Y_{n,1}), \ldots, f^{\ell}(t_{n-1}+hc_s,Y_{n,s})), \quad
\ell = 1, \ldots,D.
\]
That function $\boldsymbol{f}$ represents the parallel evaluation of $f(t_{n-1}+hc_1,Y_{n,1}), \ldots, f(t_{n-1}+hc_s,Y_{n,s})$.
For instance, if $D=2$, $s=2$, and $f(t,[y^1, y^2]) = [2\, y^1*y^2, t*(y^1+y^2)]$, then $\boldsymbol{f}(t + h \boldsymbol{c}, \boldsymbol{Y_n})$ returns the array 
\[
\left[
\begin{matrix}
2\, \boldsymbol{Y_n}^1 \circ \boldsymbol{Y_n}^2\\
(t+ h \boldsymbol{c}) \circ (\boldsymbol{Y_n}^1+\boldsymbol{Y_n}^2)
\end{matrix}
\right]
= 
\left[
\begin{matrix}
(2\, Y_{n,1}^1*Y_{n,1}^2,Y_{n,2}^1*Y_{n,2}^2)\\
((t+hc_1)*(Y_{n,1}^1+Y_{n,1}^2),(t+hc_2)*(Y_{n,2}^1+Y_{n,2}^2))
\end{matrix}
\right]
\]
Note how operations on 
$s$-vectors are performed component-wise, enabling parallel evaluation across all stages.

With that notation, (\ref{eq:Yi2}) can be written as follows, 
\begin{equation}
\label{eq:Yi3}
\begin{split}
\boldsymbol{F_n} &= \boldsymbol{f}(t_{n-1}+h \, \boldsymbol{c},\boldsymbol{Y_{n}}), \\
\boldsymbol{L_n^{\ell}} &= h \, ( \boldsymbol{b} \circ  \boldsymbol{F_n^\ell}), \quad \ell = 1,\ldots,D, \\
\boldsymbol{Y_n^{\ell}} &= y_{n-1}^{\ell} + \sum_{i=1}^{s} \boldsymbol{\mu_i}\,  L_{n,i}^{\ell}, \quad \ell = 1,\ldots,D,
\end{split}
\end{equation}
while (\ref{eq:yn2}) can be rewritten as
\begin{equation}
\label{eq:yn3}
y_{n}^{\ell} = y_{n-1}^{\ell} + \mathrm{sum}(\boldsymbol{L_n^{\ell}}), \quad \ell = 1,\ldots,D.
\end{equation}

\paragraph*{Vectorized implementation of fixed-point iteration} 

Following (\ref{eq:Yni0})--(\ref{eq:Pol0}), we initialize the fixed-point iteration by approximating $\boldsymbol{Y_n^{\ell}}$ for $\ell=1,\ldots,D$ as
\begin{equation}
\label{eq:IRK_init}
\boldsymbol{Y_n^{\ell,[0]}}= y_{n-1}^{\ell} + \sum_{i=1}^s \boldsymbol{\nu_i}\, 
 L_{n-1,i}^{\ell}.
\end{equation}
Then, $\boldsymbol{Y_{n}^{[k]}} = (Y_{n}^{1,[k]}, \ldots,Y_{n}^{D,[k]})$ is computed for $k=1,2,\ldots$ as follows:
\begin{equation}
\label{eq:IRK_FP}
\begin{split}
\boldsymbol{F_n^{[k]}} &= \boldsymbol{f}(t_{n-1}+h\, \boldsymbol{c},\boldsymbol{Y_{n}^{[k-1]}}), \\
\boldsymbol{L_n^{\ell,[k]}} &= h \, (\boldsymbol{b} \circ  \boldsymbol{F_n^{\ell,[k]}}), \quad \ell = 1,\ldots,D, \\
\boldsymbol{Y_n^{\ell,[k]}} &= y_{n-1}^{\ell} + \sum_{i=1}^{s} \boldsymbol{\mu_i}\,  L_{n,i}^{\ell,[k]}, \quad \ell = 1,\ldots,D.
\end{split}
\end{equation}
This formulation naturally enables SIMD implementation: the 
$s$-vectors can be stored in SIMD registers, and operations like $\boldsymbol{b} \circ \boldsymbol{F_n^\ell}$
correspond to single SIMD instructions.

As for the alternative IRKN form (\ref{eq:Gi})--(\ref{eq:qnvn2}) of IRKGL methods for second order ODE systems,
we initialize the fixed-point iteration as
\begin{equation}
\label{eq:IRK_init2}
\begin{split}
\boldsymbol{Y_n^{\ell,[0]}} &= y_{n-1}^{\ell} + \sum_{i=1}^s \boldsymbol{\nu_i}\,  L_{n-1,i}^{\ell}, \quad \ell=d+1,\ldots,2d,\\
\boldsymbol{L_n^{\ell,[0]}} &= h \, ( \boldsymbol{b} \circ  \boldsymbol{Y_n^{d+\ell,[0]}}), \quad \ell = 1,\ldots,d, \\
\boldsymbol{Y_n^{\ell,[0]}} &= y_{n-1}^{\ell} + \sum_{i=1}^{s} \boldsymbol{\mu_i}  L_{n,i}^{\ell,[0]}, \quad \ell = 1,\ldots,d.
\end{split}
\end{equation}
Then, for $k=1,2,\ldots$,
\begin{equation}
\label{eq:IRK_FP2}
\begin{split}
\boldsymbol{F_n^{[k]}} &= \boldsymbol{f}(t_{n-1}+h\, \boldsymbol{c},Y_n^{[k-1]}), \\
\boldsymbol{L_n^{\ell,[k]}} &= h \, ( \boldsymbol{b} \circ  \boldsymbol{F_n^{\ell,[k]}}), \quad \ell = d+1,\ldots,2d, \\
\boldsymbol{Y_n^{\ell,[k]}} &= y_{n-1}^{\ell} +
h\, \left(\boldsymbol{c} \,  y_{n-1}^{d+\ell}
+ \sum_{i=1}^{s} \boldsymbol{\eta_i}\,  L_{n,i}^{d+\ell,[k]}\right), \quad \ell = 1,\ldots,d, 
\end{split}
\end{equation}

\paragraph*{Stopping criterion}

In~\cite{Hairer2008}, it is shown that standard stopping criteria based on prescribed iteration error tolerances result in linear growth of energy error when applied to Hamiltonian systems. To address this issue, a new stopping criterion is proposed in~\cite{antonana2017}. We now present a simplified version of this criterion that is suitable for our vectorized implementation.

If $ \boldsymbol{Y_n^{\ell,[k]}} = \boldsymbol{Y_n^{\ell,[k-1]}} $ for all $ \ell \in \{1, \ldots, D\} $, the iteration should clearly be stopped; as $\boldsymbol{Y_n^{\ell,[k+1]}} = \boldsymbol{Y_n^{\ell,[k]}}$ would hold for all $ \ell \in \{1, \ldots, D\} $. However, since this condition may not always be met, the stopping criterion must be supplemented with an additional condition to detect when successive approximations cease to improve:
Define 
\begin{equation}
\label{eq:Delta_n}
\Delta_n^{\ell,[k]} = \| \boldsymbol{Y_{n}^{\ell,[k]}} - \boldsymbol{Y_{n}^{\ell,[k-1]}} \|_{\infty}. 
\end{equation}
The iteration will stop after the $k$th iteration if the following condition is met for all $\ell \in \{1, \ldots, D\}$:
\begin{equation}
\label{eq:stopping}
\Delta_{n}^{\ell,[k]} = 0 \quad \mbox{or} \quad  \min \left( \Delta_n^{\ell,[1]}, \ldots, \Delta_n^{\ell,[k-2]} \right) \leq \min \left( \Delta_n^{\ell,[k-1]}, \Delta_n^{\ell,[k]} \right).
\end{equation}
For second order ODE systems, condition (\ref{eq:stopping}) is checked  for all $\ell \in \{1, \ldots, d\}$.
 
 \section{Implementation aspects of IRKGL16}
 
\subsection{SIMD-paralellization}

Single Instruction, Multiple Data (SIMD) is a parallelization technique supported by modern CPU cores \cite{Reinders2016,Cardoso2017}. SIMD enables the execution of vectorized instructions that apply a single operation simultaneously to multiple data elements, thereby accelerating performance in computationally intensive tasks. Modern CPUs contain specialized registers known as {\em short vectors}, typically 256 bits (holding four 64-bit double-precision floating-point numbers) or 512 bits (holding eight 64-bit double-precision floating-point numbers) in size.

In a standard SIMD operation, two input vectors are processed element-wise, applying the same operation to each pair of corresponding elements and producing an output vector. By operating on multiple elements in parallel, SIMD significantly improves performance.

\subsection{Julia language and the package \texttt{SIMD.jl}}

Our SIMD-vectorized implementation of IRKGL16 is written in Julia language
 and relies on the package \texttt{SIMD.jl}~\cite{SIMD}.
Julia is a high-level dynamic language  that allows programmers to write clear, high-level, generic and abstract code  resembling  mathematical formulas, yet produces fast, low-level machine code that has traditionally only been generated by static languages  \cite{Bezanson2017, Bezanson2018}. 
By default, Julia employs just-in-time (JIT) compilation, generating LLVM intermediate code that the LLVM compiler framework then translates into optimized machine code. %specific to the target platform.

In Julia, there are several ways to explicitly apply SIMD vectorization, and the package \texttt{SIMD.jl} provides a convenient solution for this purpose. That package introduces the parameterized vector type \verb|Vec{s, T}|, representing vectors of $s$ elements of type \verb|T|. In our vectorized implementation of $s$-stage IRKGL methods (where $s = 8$ in IRKGL16), we use vectors of type \verb|Vec{s, Float64}|, meaning $s$-vectors with elements of type \texttt{Float64} (64-bit double-precision floating-point numbers), or alternatively, vectors of type \verb|Vec{s, Float32}|. The standard arithmetic and logical operations are designed to be applied element-wise in parallel, producing a SIMD vector as the result. The package \texttt{SIMD.jl} generates LLVM code that defines vectors of $s$ elements, which the LLVM compiler then translates into optimized instructions that utilize the SIMD registers available on the target platform, while also enabling other low-level optimizations (e.g., instruction-level parallelism) afforded by the \verb|Vec{s,T}| abstraction.

\subsection{VecArray and IRKGL16}

The vectorized form of fixed-point iteration (\ref{eq:IRK_init})--(\ref{eq:IRK_FP})  (resp., (\ref{eq:IRK_init2})--(\ref{eq:IRK_FP2})) for first-order problems (\ref{eq:ivp}) (resp., second-order problems (\ref{eq:ivp2})) is well-suited for SIMD implementation.  SIMD-vectorized implementation is particularly effective for IRKGL16, where $s = 8$, as this choice maximizes efficient utilization of SIMD registers on modern hardware;
 in this case, all additions and multiplications in (\ref{eq:IRK_init})--(\ref{eq:IRK_FP}) and (\ref{eq:IRK_init2})--(\ref{eq:IRK_FP2}), as well as the computations in (\ref{eq:Delta_n}), can be implemented  efficiently using SIMD operations with short vectors containing eight 64-bit double-precision floating-point numbers.  
 
In addition, the evaluation of $F(t_{n-1}+ \boldsymbol{c}\, h, \boldsymbol{Y_n^{[k-1]}})$ in (\ref{eq:IRK_FP}) and (\ref{eq:IRK_FP2}) (which is equivalent to the evaluation of $f(t_{n-1}+c_i h, Y_{n,i}^{[k-1]})$ for $i=1,\ldots,s$) also benefits from SIMD parallelization:
In IRKGL16,  parallelization remains transparent to the user, assuming the right-hand side of the system is defined in terms of arithmetic operations and elementary functions compatible with \texttt{SIMD.jl}. Such seamless SIMD parallelization 
is handled automatically with the help of  \texttt{SIMD.jl} and an additional parametrized abstract type, \verb|VecArray{s,T}|, that we have implemented to efficiently represent and handle abstract arrays of elements of type \verb|Vec{s,T}| (where \verb|T| is either \verb|Float64| or \verb|Float32|). The $s$-vectors $\boldsymbol{b,c,\mu_i,\eta_i,\nu_i}$, $\boldsymbol{Y_n^\ell}$, $\boldsymbol{F_n^\ell}$, and $\boldsymbol{L_n^\ell}$ considered in Subsection~\ref{ss:vec_IRK} are represented in IRKGL16 as vectors of type \verb|Vec{s,T}|, while $\boldsymbol{Y_n = (Y_{n}^1, \ldots,Y_{n}^D)}$, $\boldsymbol{F_n = (F_{n}^1, \ldots,F_{n}^D)}$, and 
 $\boldsymbol{L_n = (L_{n}^1, \ldots,L_{n}^D)}$, and also $\boldsymbol{\mu=(\mu_1,\ldots,\mu_s)}$, $\boldsymbol{\eta=(\eta_1,\ldots,\eta_s)}$, and  $\boldsymbol{\nu=(\nu_1,\ldots,\nu_s)}$, are represented as one-dimensional abstract arrays of type \verb|VecArray{s,T}|.
 
  Let us consider for instance the \verb|VecArray| object  \verb|Yn| representing $\boldsymbol{Y_n}$:
 \begin{itemize} 
 \item The syntax \verb|Yn[j]| (or equivalently \verb|getindex(Yn,j)|) can be used to get the vector of type \verb|Vec{s,T}| representing $\boldsymbol{Y_n^\ell}$.  
  \item The syntax \verb|Yn[j] = V| (or equivalently \verb|setindex!(Yn,j,V)|) can be used to update the $\ell$th vector of type \verb|Vec{s,T}| representing $\boldsymbol{Y_n^\ell}$.  Here, we are assuming that \verb|V| is a vector of type \verb|Vec{s,T}|.
\item The data of  \verb|Yn| are internally stored as a 2-dimensional array of type \verb|Array{T,2}| having $s$ rows and $D$ columns.
\item The $s$ components of  \verb|Yn[j]|  are stored in the $j$th column of that 2-dimensional array.
\end{itemize}

 Listing~\ref{code:ex1} shows the section of the IRKGL16 code where a single vectorized fixed-point iteration (\ref{eq:IRK_FP}) is implemented.

\begin{minipage}{\hsize}
%\lstset{frame=single, framexleftmargin=-1pt, framexrightmargin=-17pt, framesep=12pt, linewidth=0.98\textwidth, language=julia}%

\begin{lstlisting} [caption=Julia code for vectorized iteration in IRKGL16, label={code:ex1}]
f_ODE!(Fn,Yn,parms,tn+h*c) 

for j in 1:D
     Fnj=Fn[j]
     Lnj=h*(b*Fnj)
     dYnj=mu[1]*Lnj[1]
     for i in 2:s
         dYnj=dYnj+mu[i]*Lnj[i]
     end
     Yn[j]=yn[j]+dYnj
end


\end{lstlisting}

\end{minipage}
  
\medskip

The first line of this code requires some explanation, while the remaining lines in Listing~\ref{code:ex1} directly translate the last two lines of (\ref{eq:IRK_FP}). The function \verb|f_ODE!(dy, y, parms, t)| is a user-defined generic function that takes as inputs the array \verb|y|, a vector \verb|parms| of constant parameters, and the time \verb|t|. It computes the $ D $ components of the right-hand side $ f(t, y) $ of the ODE system, storing the result in the pre-existing array \verb|dy|. For illustration, Listing~\ref{code:ex2} provides an implementation of \verb|f_ODE!(dy, y, parms, t)| for the following non-autonomous Hamiltonian ODE system: 
\begin{align}
\begin{split}
\frac{dq_1}{dt}&=p_1, \quad 
\frac{dq_2}{dt}=p_2, \\
\frac{dp_1}{dt}&=-q_1-2 (\lambda + \xi \sin(t)) q_1q_2, \\
\frac{dp_2}{dt}&=-q_2-(\lambda + \xi \sin(t)) (q_1^2-q_2^2).
\end{split}
\end{align}
Note that this system reduces to the classic Hénon-Heiles system when $\xi=0$.
\begin{minipage}{\hsize}
	%\lstset{frame=single,framexleftmargin=-1pt,framexrightmargin=-17pt,framesep=12pt,linewidth=0.98\textwidth,language=julia}% 
	%
\begin{lstlisting} [caption=User-defined generic function for the Hénon-Heiles problem, label={code:ex2}]
function f_ODE!(dy,y,parms,t)

	# q1=y[1]; q2=y[2]; p1=y[3]; p2=y[4]
	lambda=parms[1]
	xi=parms[2]
	aux=lambda+xi*sin(t)
	dy[1]=y[3]
	dy[2]=y[4]
	dy[3]=-y[1]-2*aux*y[1]*y[2]
	dy[4]=-y[2]-aux*(y[1]^2-y[2]^2)
	
	return nothing
end

\end{lstlisting}

\end{minipage}

 The generic function \verb|f_ODE!(dy, y, parms, t)| in Listing~\ref{code:ex2} supports a wide range of argument types. For example, it can be called with \verb|t| as a scalar of type \verb|T = Float64| (or alternatively \verb|T=Float32| or \verb|T=BigFloat|), and \verb|dy| and \verb|y| as standard Julia arrays containing elements of type \verb|T|. Additionally, the function accepts other argument types, provided that the operations within the \verb|f_ODE!| definition are valid for them. Specifically, the function call \verb|f_ODE!(Fn, Yn, parms, tn + h*c)| in Listing~\ref{code:ex1}, where \verb|Fn| and \verb|Yn| are abstract arrays of type \verb|VecArray{s, T}| (with \verb|T = Float64| or \verb|T = Float32|) and \verb|tn + h*c| is of type \verb|Vec{s, T}|, is equivalent to the first line of (\ref{eq:IRK_FP}), representing the parallel evaluation of $f(t_{n-1}+c_i h, Y_{n,i}^{[k-1]})$ for $i=1,\ldots,s$.
This design enables parallel evaluations of the ODE system’s right-hand side in a manner that is transparent to the user, provided that the right-hand side \(f(t,y)\) is implemented as a generic function 
\verb|f_ODE!(dy, y, parms, t)| 
without branching and using only arithmetic operations and elementary functions compatible with  \texttt{SIMD.jl}. 
All test problems considered in Section~\ref{sec5} satisfy this requirement.

However, not all ODE systems~(\ref{eq:ivp}) admit a generic implementation of \verb|f_ODE!(dy, y, parms, t)| that supports seamless \(s\)-fold parallel evaluation of the system’s right-hand side. In particular, this may fail when \(f(t,y)\) contains unavoidable branching, which breaks the uniform instruction flow required for SIMD processing. In such cases, the right-hand side must be evaluated sequentially for each stage vector \(Y_i\), although the remaining vectorized computations in (\ref{eq:IRK_init})--(\ref{eq:IRK_FP}) (for first-order ODEs),  (\ref{eq:IRK_init2})--(\ref{eq:IRK_FP2}) (for second-order ODEs) and (\ref{eq:Delta_n}) remain applicable.

In the current implementation of \textsf{IRKGL16}, simultaneous SIMD evaluation of the user-supplied function $f$ is performed only when the integrator option \verb|fseq=false| is explicitly selected by the user. All other computations are executed using SIMD vectorization whenever the state variables are of type \verb|Float32| (single precision) or \verb|Float64| (double precision). A fully sequential fallback implementation is used for state variables of other types (e.g., \verb|BigFloat|) or whenever the \verb|simd|   option is explicitly disabled (\verb|simd=false|) by the user.

 %however, IRKGL16 can still be applied if the user provides a specialized implementation of \verb|f_ODE!(dy, y, parms, t)| for instances where \verb|dy| and \verb|y| are abstract arrays of type \verb|VecArray{s, T}| (with \verb|T = Float64| or \verb|T = Float32|) and \verb|t| is of type \verb|Vec{s, T}|.

  In some cases, such as the $N$-body example considered in Section~\ref{sec5}, it is more convenient for the user to organize the state variables in a higher-dimensional array rather than a one-dimensional array. To accommodate this, we have parameterized our abstract array type with an additional parameter \verb|dim|, as \verb|VecArray{s, T, dim}|. For one-dimensional arrays of vectors of type \verb|Vec{s, T}|, we set \verb|dim = 2|, with the underlying data stored in a two-dimensional standard array of type \verb|Array{T, 2}|. More generally, \verb|VecArray{s, T, dim}| refers to $(\texttt{dim} - 1)$-dimensional abstract arrays of vectors of type \verb|Vec{s, T}|, with the data stored in a \verb|dim|-dimensional array of type \verb|Array{T, dim}|. Our implementations of the \verb|getindex| and \verb|setindex!| functions work as expected in this more general case, allowing users to apply IRKGL16 with a generic implementation of the function \verb|f_ODE!(dy, y, parms, t)| when the state variables are stored in a higher-dimensional array.

  \section{State-of-the-art explicit symplectic integrators}

  Splitting and composition techniques are powerful tools for constructing practical symplectic integrators in many areas of application. Although the resulting methods are typically tailored to the specific problem and lack the general applicability of approaches such as IRKGL schemes, they can be highly efficient when applicable.
  \cite{Sanz1994,Blanes2024,Hairer2006}.
  
  In the context of Hamiltonian systems that can be decomposed into two (or more) exactly solvable components, explicit symplectic integrators can be constructed via operator splitting techniques. A widely used approach is the \emph{second-order Strang splitting method}~\cite{Strang1968}, which yields a symmetric and symplectic scheme by composing the $t$-flow maps $\varphi^A_{t}$ and $\varphi^B_{t}$ associated with the individual sub-Hamiltonians. Specifically,  the Strang splitting scheme advances the solution over a time step $h$ via the composition
\begin{equation}
\label{eq:Strang}
\phi_h = \varphi^A_{h/2} \circ \varphi^B_{h} \circ \varphi^A_{h/2}.
\end{equation}
This method is second-order accurate, time-reversible, and symplectic.  The resulting integrator is explicit when the Hamiltonian is written as $H = H_A + H_B$, and the $t$-flows $\varphi^A_t$ and $\varphi^B_t$ generated by $H_A$ and $H_B$, respectively, are available in closed form.

For instance, consider a separable Hamiltonian of the form~\eqref{eq:sHam}, where $H_A(q,p) = \tfrac{1}{2} p^T M^{-1} p$ and $H_B(q,p) = U(q)$. Via the linear change of variables $p = M v$, the exact flows are given by 
\begin{equation}
\label{eq:flowAB}
\varphi^A_t(q,v) = (q + t \, v), \quad 
\varphi^B_t(q,v) = (q,\, v - t \, g(q)), 
\end{equation}
where $g(q) =  M^{-1} \nabla U(q)$.

The method can be extended to systems where the Hamiltonian splits into more than two exactly solvable components. For instance, if $H = H_A + H_B + H_C$ then a symmetric second-order integrator can be obtained by composing the flows of each part as
\begin{equation}
\label{eq:genStrang}
\phi_h = \varphi^A_{h/2} \circ \varphi^B_{h/2} \circ \varphi^C_{h} \circ \varphi^B_{h/2} \circ \varphi^A_{h/2}.
\end{equation}
This generalized Strang-type splitting remains symplectic and time-symmetric, provided that each subflow $\varphi^A_t$, $\varphi^B_t$, and $\varphi^C_t$ is symplectic. Such multi-part splittings arise naturally in applications like charged particle dynamics, Lie-Poisson systems, and systems with coupled fast and slow Hamiltonian components~\cite{Blanes2024, McLachlanQuispel2002}.

    Extensions to higher-order accuracy can be obtained by symmetric compositions of $s$ basic Strang stages ~\cite{Blanes2024,Hairer2006} of the form,
  \begin{align}
  \label{eq:comp_method}
  &\psi_{h}= %\phi_{\gamma_{1}h} \circ \phi_{\gamma_{2}h} \circ \cdots \circ  \phi_{\gamma_{s-1}h} \circ 
  \phi_{\gamma_{s}h} \circ \phi_{\gamma_{s-1}h} \circ \dots \circ \phi_{\gamma_{2}h} \circ \phi_{\gamma_{1}h}
  \end{align}
  with $\gamma_{s+1-j}=\gamma_{j}$, $1\leq j \leq s$.
 
%     In the case of descomposition into two solvable Hamiltonians, we have a particular explicit symplectic integrator knows as Runge-Kutta-Nystr\"{o}m methods where the separable Hamiltonian is the form,
    
      We consider the following symmetric compositions of time-symmetric second-order schemes of orders $r = 2, 4, 6, 8, 10$, as recommended in \cite{Blanes2024}:
\begin{itemize}
    \item[1)] Strang splitting: a second-order method ($r = 2$) with $s = 1$ stage.
    \item[2)] SUZ90: a fourth-order method ($r = 4$) with $s = 5$ stages, proposed by Suzuki (1990) \cite{Suzuki1990}.
    \item[3)] SS05:  a sixth-order scheme with $s = 13$ stages, an eighth-order one with $s = 21$ stages, and a tenth-order one with $s = 35$ stages proposed by Sofroniou and Spalleta (2005)~\cite{Sofroniou2005}.
\end{itemize}

  Any composition of the form \eqref{eq:comp_method} with \eqref{eq:Strang} can be rewritten as
  \begin{equation}
  \label{eq:high_order_splitting}
\varphi_{a_{s+1} h}^{A} \varphi_{b_s h}^{B} \circ  \varphi_{a_s h}^{A} \circ \dots \circ  \varphi_{b_1 h}^{B} \circ  \varphi_{a_1 h}^{[A]}
\end{equation}
with coefficients
\begin{equation*}
a_1 = \frac{\gamma_1}{2},  \quad a_j = \frac{\gamma_{j-1}+\gamma_{j}}{2}, \quad
b_j = \gamma_j. 
\end{equation*}
Clearly,  \eqref{eq:high_order_splitting} cannot be written in the form \eqref{eq:comp_method} for arbitrary sequences of coefficients $a_j,b_j$. Such more general family of methods are of particular interest
in the case of Hamiltonians of the form \eqref{eq:sHam}, for which efficient high-order schemes have been constructed. They are often referred as splitting methods of Runge-Kutta-Nystr\"{o}m (RKN) type, as they can be rewritten as explicit symplectic RKN methods applied to the system of second order ODEs \eqref{eq:ivp2}~\cite{Sanz1994}.
   We consider the following splitting integrators of RKN type of order $r=6$ and $r=8$ recommended in \cite{Blanes2024}:
  \begin{itemize}
  	\item[1)] $\mathrm{S_{ABA}\text{-}BM02}$: a sixth-order method with $s=14$ stages method proposed in Blanes and Moan(2002) \cite{Blanes2002}.
  	\item[2)] $\mathrm{S_{ABA}\text{-}BCE22}$: optimized eighth-order with $s=19$ stages proposed in Blanes, Casas and Escorihuela (2022) \cite{Blanes2022b}.
  \end{itemize}
    
 % \subsection{Symplectic-like extended phase integrators}
  
 % Recently, several variants of explicit integrators for the long-term integration of Hamiltonian systems without a specific separable structure have been proposed using the extended phase space technique~\cite{Pihajoki[15], Tao[16], Ohsawa[14], Luo[7]}. The symplectic-like extended phase integrators with midpoint projection introduced in~\cite{Luo[7]} have demonstrated near-conservation of energy and quadratic invariants in astrophysical applications, similar to traditional symplectic methods—hence the term "symplectic-like." In~\cite{McLachlan2023}, it is shown that the integrators proposed in~\cite{Luo[7]} are, in fact, explicit Runge-Kutta methods and that they are pseudo-symplectic and pseudo-symmetric of order $2r + 1$, provided that they have  order of accuracy $r$.

\section{Numerical experiments}
\label{sec5}

We compare the performance of IRKGL16 with several state-of-the-art explicit symplectic integrators for the high-precision integration of non-stiff Hamiltonian systems using constant step size. A range of representative examples is considered to assess both accuracy and efficiency.

It is well known~\cite{Hairer2006} that, for long-term integrations, it is advisable to apply some form of improved summation (e.g., Kahan’s compensated summation~\cite{Kahan1965}) in (\ref{eq:yn2}) to prevent excessive accumulation of roundoff errors. The idea is to increase the accuracy of the sum by representing each component of $y_n$ (for $n \geq 0$) as the sum of two floating-point numbers while representing each component of the update $\sum_{i=1}^{n} L_{n,i}$ (expected to be of relatively smaller size) as a floating-point number.  In IRKGL16, we employ Julia’s \verb|Base.TwicePrecision| for this purpose, which provides greater accuracy (though at higher computational cost) than Kahan’s compensated summation.

Some of the explicit symplectic integrators used in our comparison are available through the Julia \texttt{DifferentialEquations.jl} suite~\cite{Rackauckas2017}, but no compensated summation is applied there. For consistency and fair benchmarking, we implemented in Julia all the reference integrators described in the previous section  using (optionally) Kahan’s compensated summation. We verified that, for the problems considered, our implementations of the splitting methods are at least as efficient as (and often more efficient than) the implementations in \texttt{DifferentialEquations.jl} when both are used without compensated summation.

%In our numerical experiments, compensated summation was applied for the splitting methods in the two second-order ODE systems, whereas no compensated summation was used for the first-order ODE example to avoid additional complexity.

The source code for the integrators, as well as all Jupyter notebooks used to reproduce the numerical results, is publicly available in the repository associated with this article~\cite{supplementary_code}.

  % and show, that IRKGL16 SIMD-vectorized implementation (IRKGL16) significantly enhance the performance of the sequential IRKGL16 implementation (IRKGL16-SEQ). 

All numerical experiments were conducted on an 11th Gen Intel Core i7-11850H processor (2.5 GHz, 16 threads) equipped with 512-bit SIMD registers. We used the long-term support (LTS) stable release of Julia v1.10.10 (June 27, 2025) for all tests. 
%To optimize SIMD performance, Julia was launched from the terminal with bounds checking disabled using the flag \texttt{--check-bounds=no}.

\subsection{Test problems}

We have made comparisons for a Hamiltonian system  of first-order differential equations split into three solvable parts, and two Hamiltonian systems of second-order differential equations with Hamiltonian (\ref{eq:sHam}) separable in two solvable parts.

\subsubsection{System of first-order differential equations}

We consider a Hamiltonian system describing the motion of charged particles in the vicinity of a Schwarzschild black hole under the influence of an external magnetic field, as studied in~\cite{Naying2022note}:
\begin{align}
\begin{split}
H(r, \theta, p_r, p_{\theta}) &= \frac{1}{2} \left(1 - \frac{2}{r} \right) p_r^2 
- \frac{1}{2} \left(1 - \frac{2}{r} \right)^{-1} E^2 
+ \frac{p_{\theta}^2}{2r^2} \\
&\quad + \frac{1}{2r^2 \sin^2 \theta} \left(L - \frac{\beta}{2} r^2 \sin^2 \theta \right)^2.
\end{split}
\label{eq:Ham_SBH}
\end{align}

This Hamiltonian admits multiple splittings into analytically solvable components---typically into three, four, or five parts. Following the findings of~\cite{Naying2022note}, we adopt the three-part splitting
\[
H = H_A + H_B + H_C,
\]
which was shown to offer the best numerical performance across both regular and chaotic trajectories. The individual components are defined as follows:
\begin{align*}
H_A(r, \theta) &= \frac{1}{2r^2 \sin^2 \theta} \left(L - \frac{\beta}{2} r^2 \sin^2 \theta \right)^2 - \frac{1}{2} \left(1 - \frac{2}{r} \right)^{-1} E^2, \\
H_B(r, p_r, p_\theta) &= \frac{1}{2} \left(p_r^2 + \frac{p_{\theta}^2}{r^2} \right), \\
H_C(r, p_r) &= -\frac{1}{r} p_r^2.
\end{align*}

Although~\cite{Naying2022note} provides analytical expressions for the $t$-flow maps $\varphi_t^A$, $\varphi_t^B$, and $\varphi_t^C$, their direct implementation is computationally inefficient. To enable a fair and representative comparison with our implicit symplectic integrator (IRKGL16), we have improved the evaluation of these $t$-flows, significantly enhancing the performance of the explicit method. Unlike the analytical expressions provided in~\cite{Naying2022note} for the $t$-flows $\varphi_t^A$, $\varphi_t^B$, and $\varphi_t^C$, our implementation avoids the use of trigonometric and inverse trigonometric functions entirely, resulting in improved computational efficiency. The derivation and implementation of the optimized  flows are provided in the documentation available in the software repository~\cite{supplementary_code}  associated with this article.

As in~\cite{Naying2022note}, we set the parameters to $E = 0.995$, $L = 4.6$, and $\beta = 8.9 \times 10^{-4}$. The initial conditions are $\theta = \pi/2$ and $p_r = 0$, with $r = 11$ chosen so that the corresponding orbit is regular (i.e., lies on a closed curve). The initial value of $p_\theta > 0$ is determined by imposing the energy condition $H = -\tfrac{1}{2}$. The integration is performed over the interval $[0,10^5]$.

\subsubsection{Systems of second-order differential equations}

We consider two second-order initial value problems of the form~\eqref{eq:ivp2}, corresponding to Hamiltonian systems of the type~\eqref{eq:sHam}: the outer Solar System model~\cite{Hairer2006} and the Hénon-Heiles problem~\cite{Hairer2008}.

\paragraph*{6-body outer Solar System.}

We examine a simplified model of the outer Solar System, which includes the Sun, the four outer planets (Jupiter, Saturn, Uranus, and Neptune), and Pluto, as point masses interacting via mutual Newtonian gravity. The system is described by an 18-degree-of-freedom Hamiltonian, with $q_i, p_i \in \mathbb{R}^3$ for $i = 1, \dots, 6$, and Hamiltonian function
\begin{equation}
\label{eq:Ham2}
H(q, p) = \sum_{i=1}^N \frac{\|p_i\|^2}{2m_i} - G \sum_{1 \le i < j \le N} \frac{m_i m_j}{\|q_i - q_j\|},
\end{equation}
where $G$ is the gravitational constant, $m_i$ denotes the mass of body $i$, and $N = 6$ is the total number of bodies considered.

Initial conditions are taken from the DE430 planetary ephemerides at Julian date 2440400.5 (June 28, 1969)~\cite{Folkner2014}, and adjusted so that the barycenter of the system is at rest. The integration is carried out over a time interval of $10^7$ days.

\paragraph*{Hénon-Heiles system.}

We also consider the classical Hénon-Heiles Hamiltonian:
\begin{equation}
H(q, p) = \frac{1}{2}(p_1^2 + p_2^2) + \frac{1}{2}(q_1^2 + q_2^2) + q_1^2 q_2 - \frac{1}{3} q_2^3.
\end{equation}

This system is integrated from initial conditions corresponding to regular (non-chaotic) motion over a long time interval of length $2\pi \times 10^4$. Specifically, we choose $q_1(0) = 0$, $q_2(0) = 0.3$, $p_2(0) = 0.2$, and determine $p_1(0) > 0$ such that $H = 1/12$ (see~\cite{Hairer2006, Hairer2008}).

\subsection{Efficiency comparison of different implementations of $s$-stage IRKGL schemes}
\label{ss:IRKGL_comp}

\begin{figure}[t]
	\centering	
	\begin{subfigure}{0.48\textwidth}
		\centering
		\includegraphics[width=\textwidth]{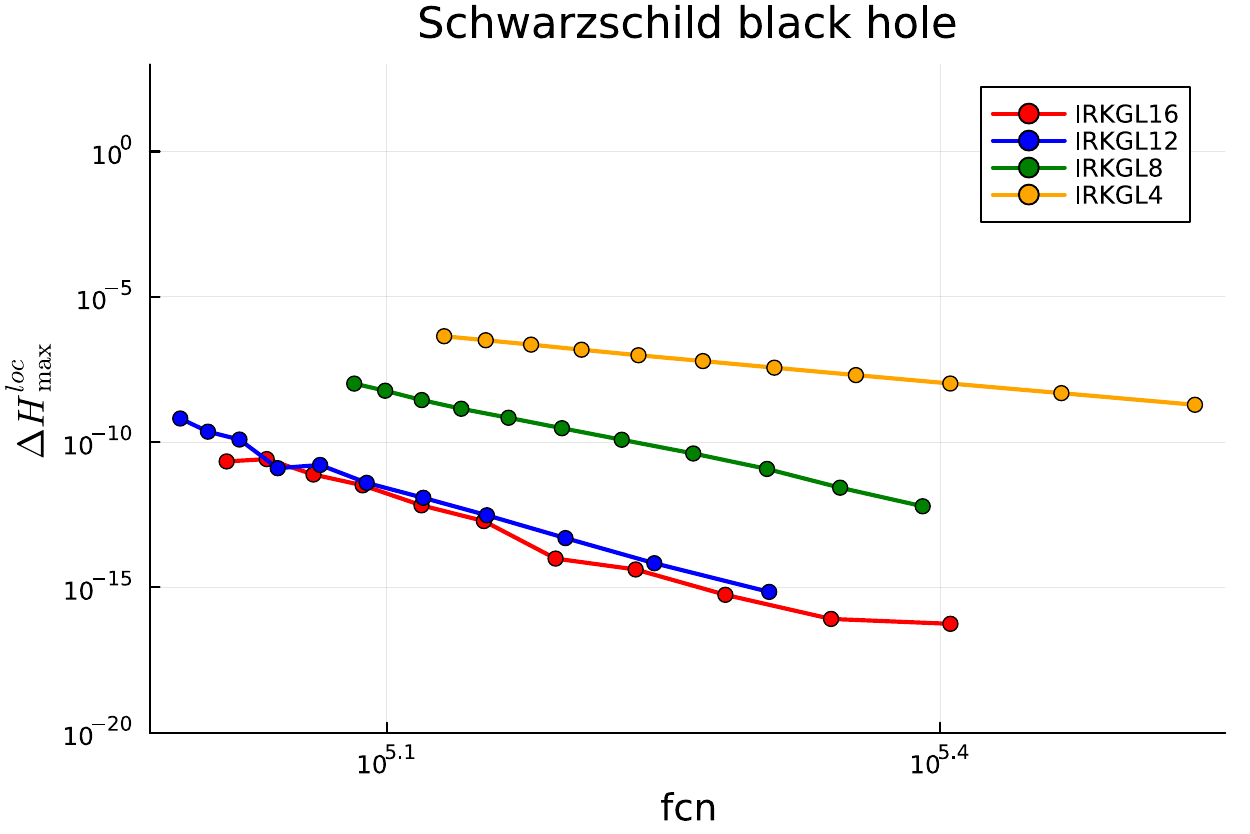}
		\caption{}
		\label{fig:ODE1_IRKGLa}
	\end{subfigure}%
	\hfill
	\begin{subfigure}{0.48\textwidth}
		\centering
		\includegraphics[width=\textwidth]{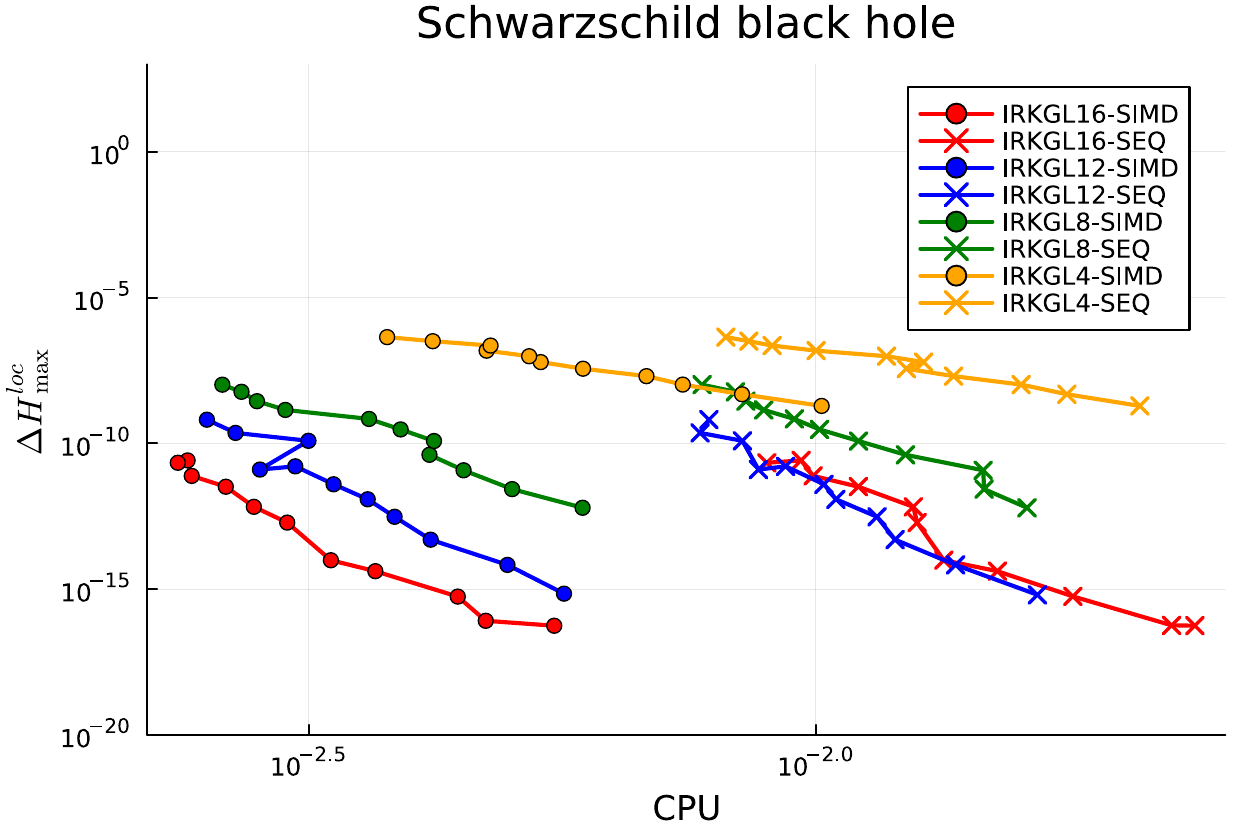}
		\caption{}
		\label{fig:ODE1_IRKGLb}
	\end{subfigure}%
	\caption{(a) Maximum local Hamiltonian error versus total function evaluations for the $s$-stage IRKGL schemes ($s = 2, 4, 6, 8$: IRKGL4, IRKGL8, IRKGL12, IRKGL16).  
(b) The same error versus CPU time for two implementations of each scheme: fully vectorized (\texttt{simd=true, fseq=false}) and fully sequential (\texttt{simd=false, fseq=true})
	  }
	\label{fig:ODE1_IRKGL}
\end{figure}

In Figure~\ref{fig:ODE1_IRKGL} we show work--precision diagrams for several implementations of $s$-stage IRKGL schemes applied to the Schwarzschild black hole problem.  

The left panel of each diagram reports the maximum local error in the Hamiltonian,
\begin{align}
\label{eq:DHlocmax}
\Delta H^{\mathrm{loc}}_{\max}
= \max_n \left( \Delta H^{\mathrm{loc}}_n \right), 
\qquad 
\Delta H^{\mathrm{loc}}_n
= \left|\frac{H(y_n)-H(y_{n-1})}{H(y_{n-1})}\right|,
\end{align}
as a function of the total number of function evaluations, for our implementations of the $s$-stage IRKGL schemes with $s=2,4,6,8$ (IRKGL4, IRKGL8, IRKGL12, and IRKGL16, respectively).  
This provides a measure of the \emph{intrinsic relative efficiency} of the methods, independent of the computational environment, under the assumptions that (i) the user-supplied function that evaluates the right-hand side of the ODE system is executed once per stage vector, in sequential order, and (ii) the total CPU time is dominated by these function evaluations.  
According to this metric, IRKGL16 is the most efficient method in the high-accuracy regime, whereas IRKGL12 performs slightly better than IRKGL16 at moderate accuracy levels.

The right panel plots the same error against actual CPU time for two versions of each $s$-stage IRKGL method: the fully vectorized implementation (\verb|simd=true, fseq=false|) and the fully sequential one (\verb|simd=false, fseq=true|).  
The relative performance of the fully sequential versions is consistent with the behavior observed in the left panels.  
In contrast, the fully vectorized implementation of IRKGL16 outperforms the lower-order schemes even at moderate accuracy.  
This improvement is due to the higher SIMD speed-up achieved by the 8-stage scheme.

\begin{figure}[t]
	\centering	
	\begin{subfigure}{0.48\textwidth}
		\centering
		\includegraphics[width=\textwidth]{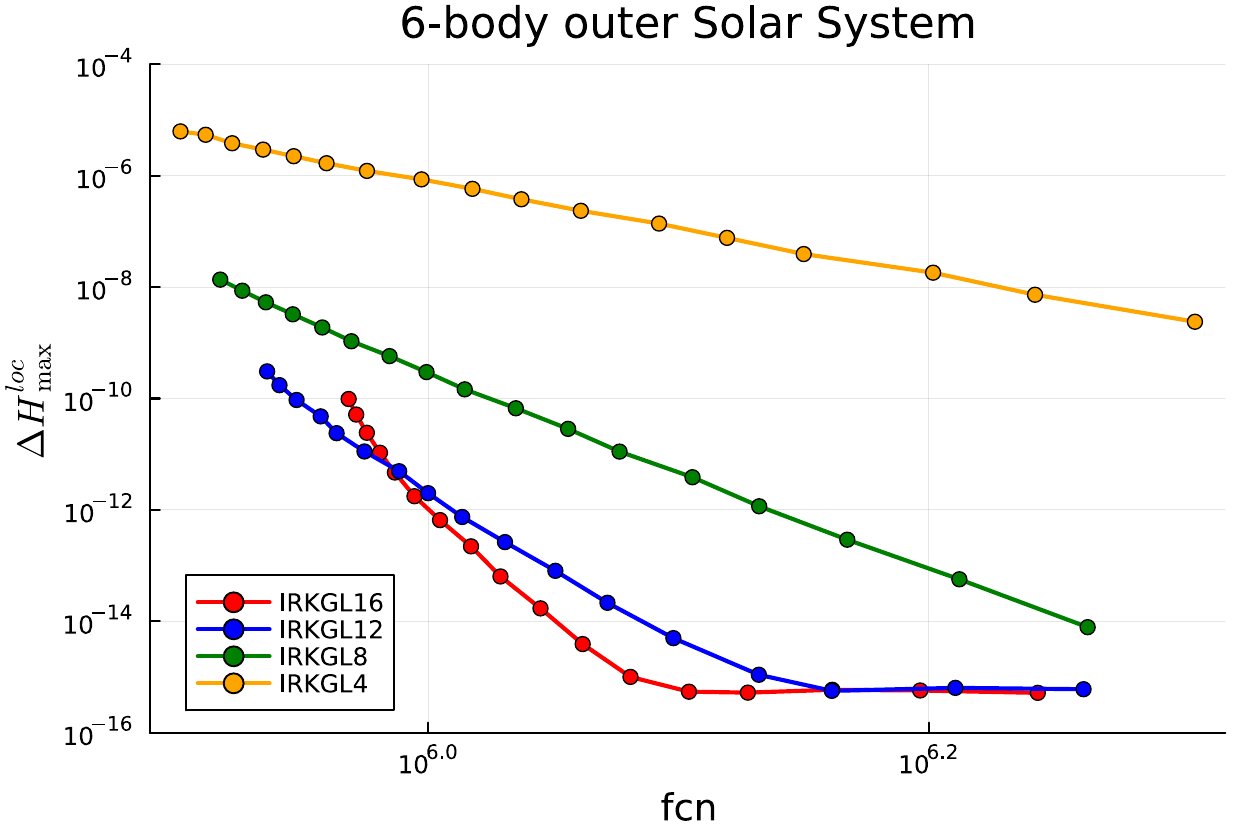}
		\caption{}
		\label{fig:ODE2_IRKGLa}
	\end{subfigure}%
	\hfill
	\begin{subfigure}{0.48\textwidth}
		\centering
		\includegraphics[width=\textwidth]{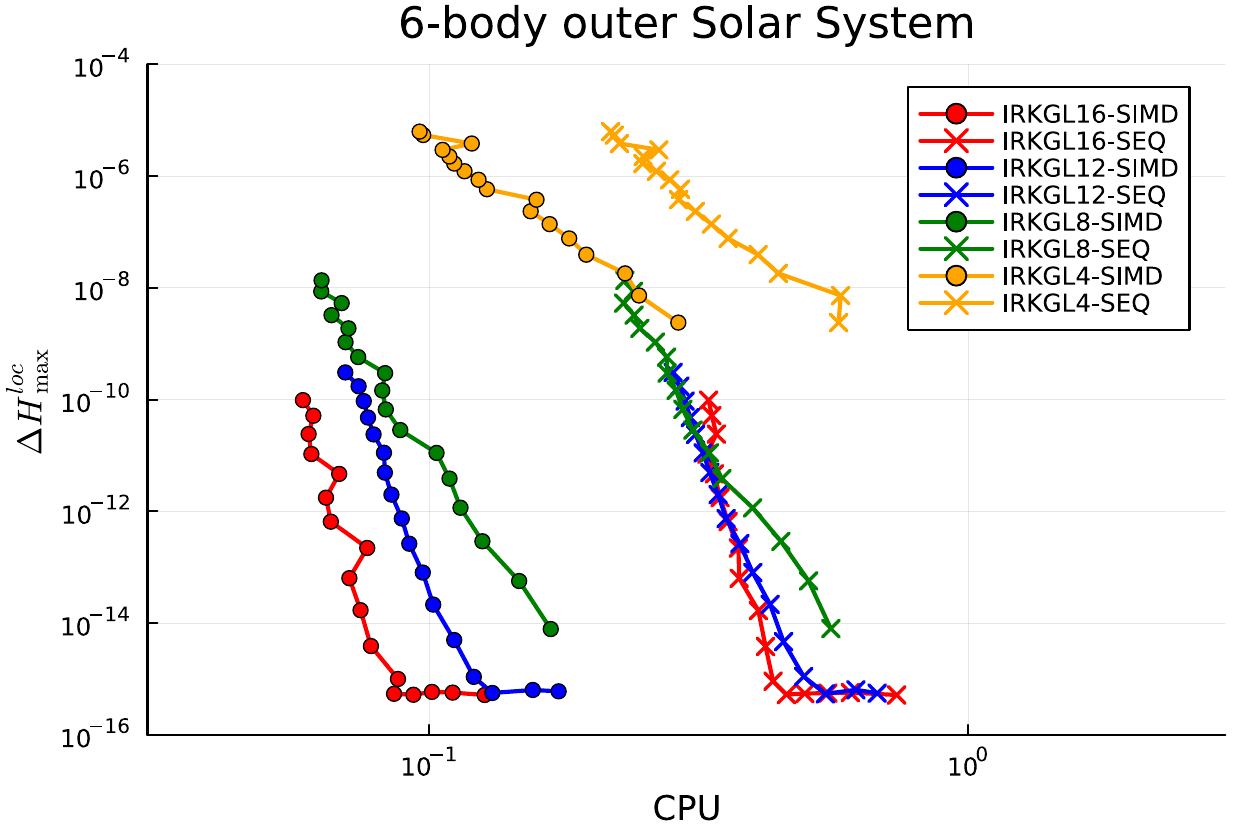}
		\caption{}
		\label{ODE2_IRKGLb}
	\end{subfigure}%
	\hfill
	\begin{subfigure}{0.48\textwidth}
		\centering
		\includegraphics[width=\textwidth]{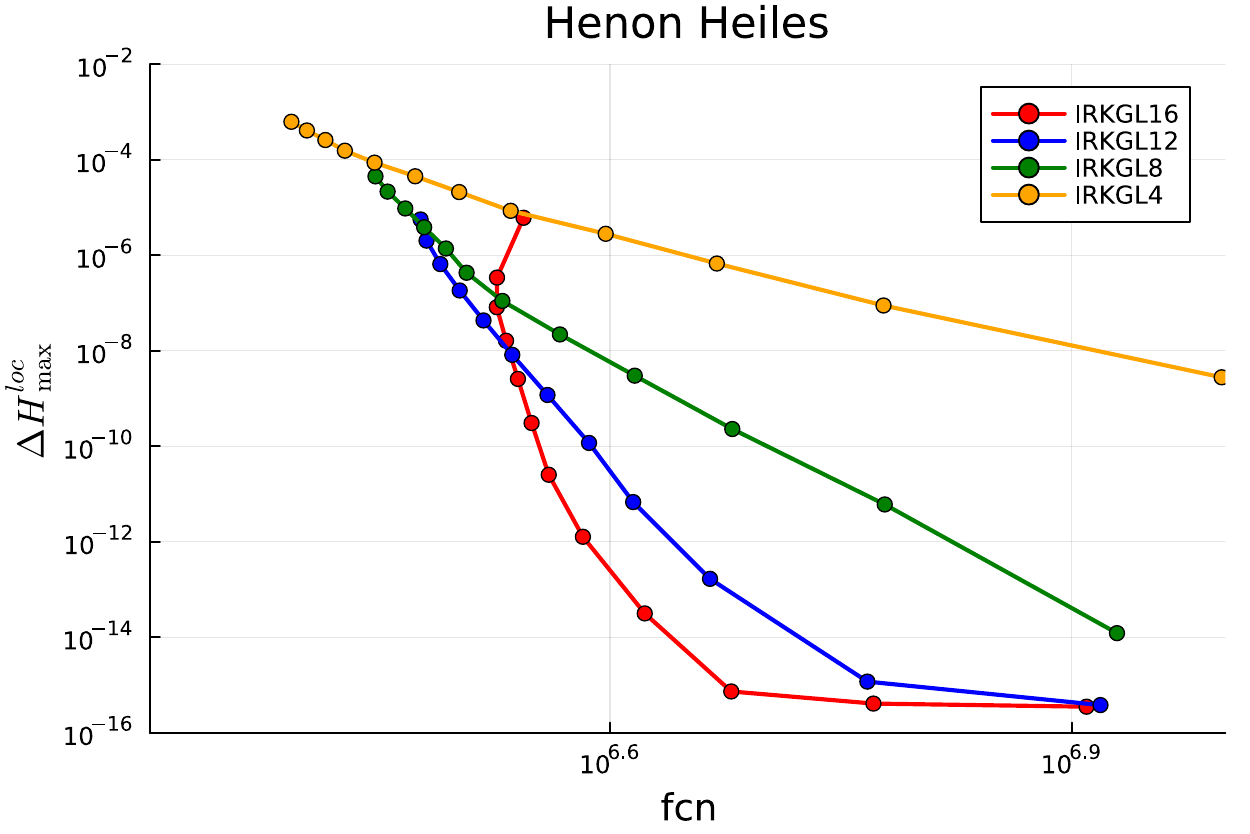}
		\caption{}
		\label{ODE2_IRKGLc}
	\end{subfigure}%
	\hfill
	\begin{subfigure}{0.48\textwidth}
		\centering
		\includegraphics[width=\textwidth]{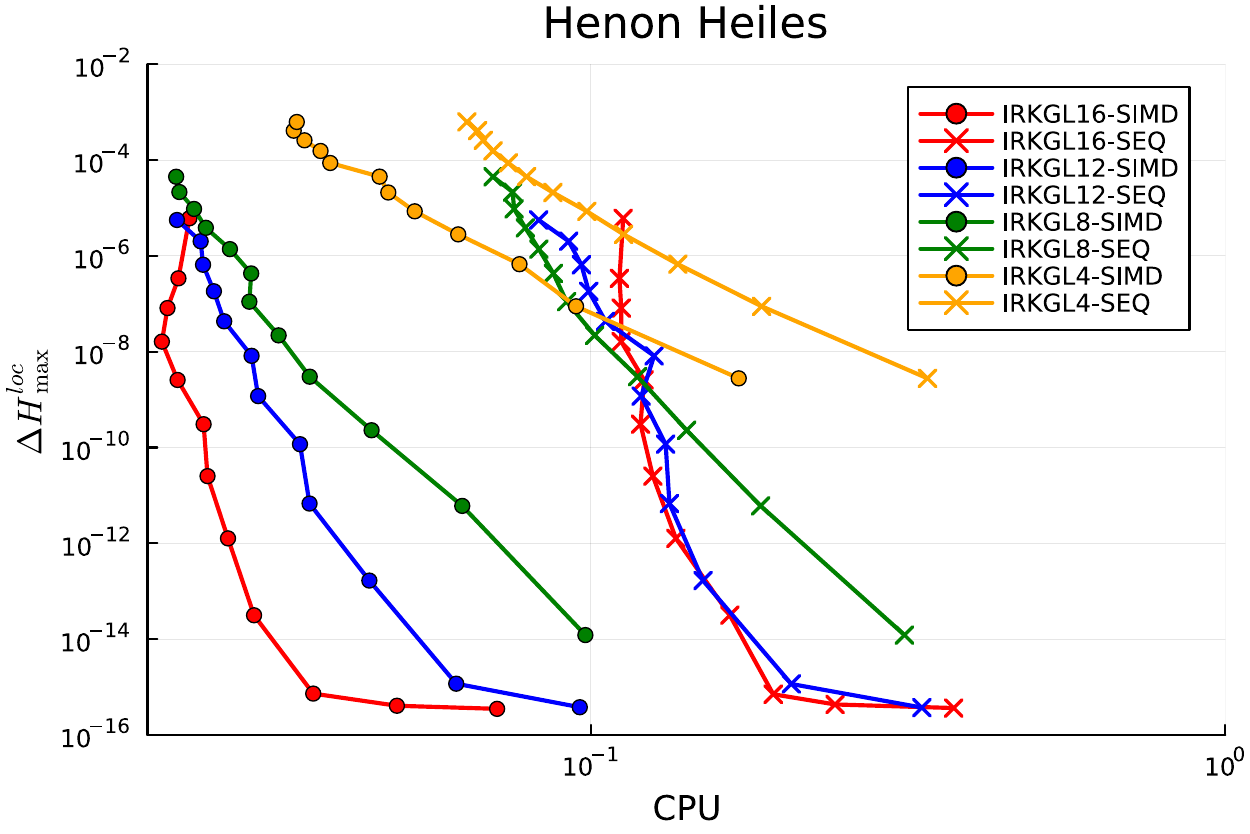}
		\caption{}
		\label{ODE2_IRKGLd}
	\end{subfigure}%
	\caption{Panels (a)–(b) correspond to the 6-body outer Solar System, and panels (c)–(d) to the Hénon–Heiles system.  
(a, c) Maximum local Hamiltonian error versus total function evaluations for the $s$-stage IRKGL schemes ($s = 2, 4, 6, 8$).  
(b, d) The same error versus CPU time for two implementations: fully vectorized and fully sequential
}
	\label{fig:ODE2_IRKGL}
\end{figure}

In Figure~\ref{fig:ODE2_IRKGL}, we present analogous work--precision diagrams for the two second-order Hamiltonian problems. The results again indicate that the fully vectorized implementation of IRKGL16 is the most efficient integrator for both moderate and high accuracy levels.

\subsection{Numerical comparisons of IRKGL16 with explicit symplectic integrators}

\subsubsection{Work--precision diagrams}

\begin{figure}[t]
	\centering	
	\begin{subfigure}{0.45\textwidth}
		\centering
		\includegraphics[width=\textwidth]{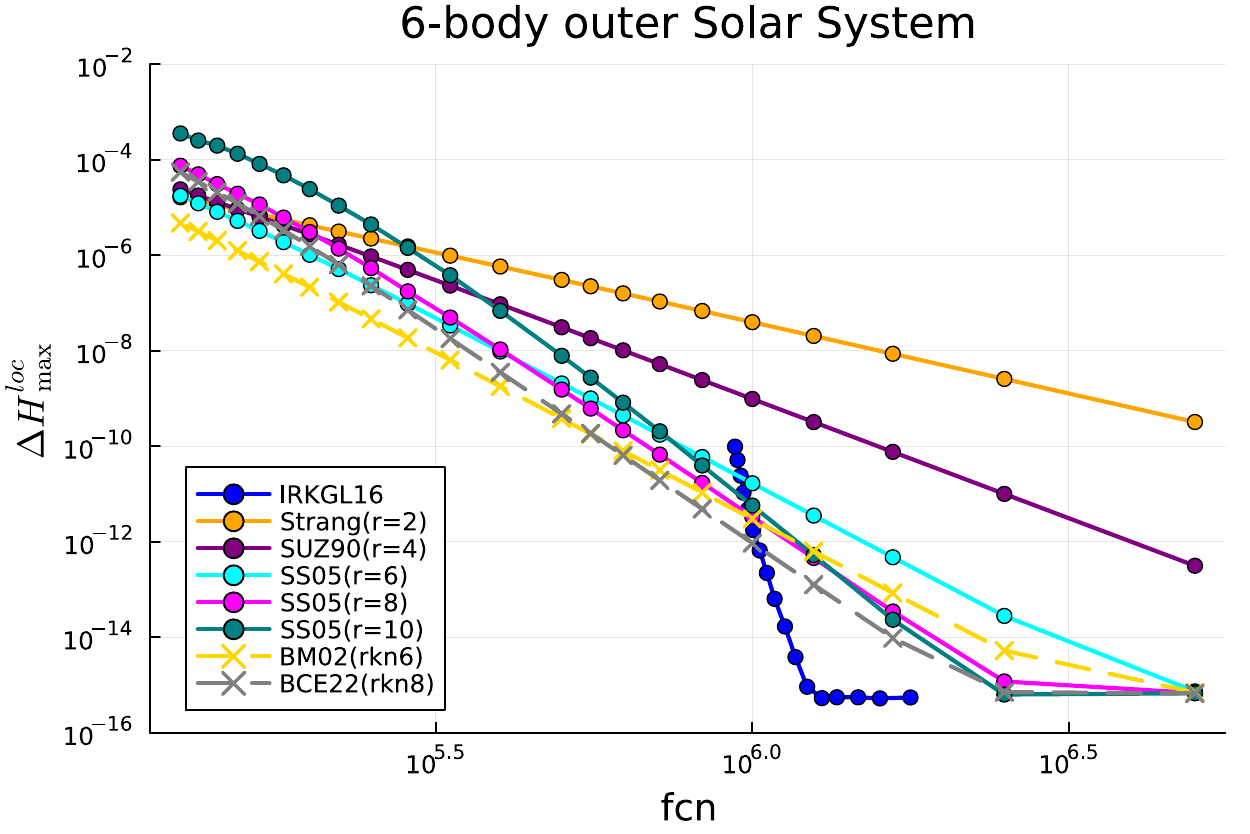}
		\caption{}
		\label{fig:Fig3a}
	\end{subfigure}%
		\hfill
	\begin{subfigure}{0.45\textwidth}
		\centering
		\includegraphics[width=\textwidth]{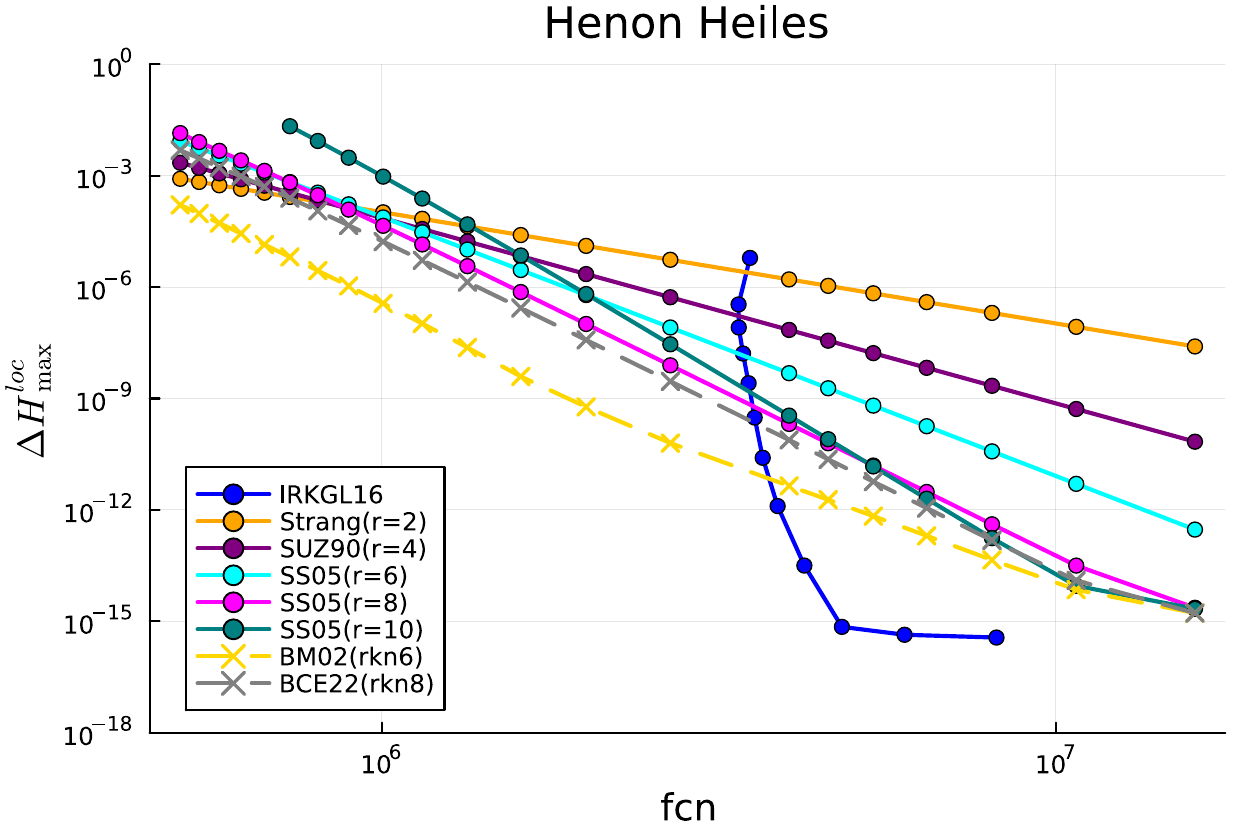} 
		\caption{}
		\label{fig:Fig3c}
	\end{subfigure}%
			\caption{Work-precision diagrams (maximum local Hamiltonian error versus number of evaluations of $g(t,u)$) for the 6-body outer Solar System and the Hénon–Heiles problem, comparing IRKGL16 with symmetric composition schemes (orders $r = 2, 4, 6, 8, 10$) and symplectic RKN methods of orders 6 and 8}	\label{fig3}
\end{figure}

We now present work--precision diagrams for the two second-order Hamiltonian problems 
considered in Subsection 5.1.2. The computational cost of one step of any splitting 
method~(\ref{eq:high_order_splitting}) applied with the flows \eqref{eq:flowAB} is 
dominated by the \(s\) applications of the \(t\)-flow 
\(\varphi^{B}_t(q,v) \mapsto (q, v + t\, g(q))\), where \(g(q)=M^{-1}\nabla U(q)\). Thus, for splitting methods, the 
relative cost of a numerical integration can essentially be measured by the number of 
calls to the function \(g(u)\). For implicit methods applied to Hamiltonian second-order 
ODEs of the form (\ref{eq:Hamode2}), the relative cost can also be quantified in terms of 
function evaluations of the right-hand side \(g(q)\), provided that it 
is sufficiently expensive to compute.

Figure~\ref{fig3} presents work--precision diagrams for the two second-order ODE problems, 
displaying the maximum local error in the Hamiltonian \eqref{eq:DHlocmax} plotted against 
the number of function evaluations for several splitting methods and IRKGL16. We consider 
the splitting integrators obtained by applying the composition schemes of the previous 
subsection with (\ref{eq:Strang}) as a basic integrator, and additionally the two RKN-type 
splitting methods of orders 6 and 8 selected in Section~4, all implemented with 
compensated summation.

From the efficiency diagrams shown in Figure~\ref{fig3}, the larger overhead inherent to implicit schemes implies that the best explicit symplectic methods will outperform sequential implementations of IRKGL16 at low and moderate accuracy levels, regardless of the programming language or other aspects of the computational environment. 
However, at higher target accuracies,  i.e., when the maximum local energy error is just above roundoff level, a sequential implementation of IRKGL16 may surpass explicit methods, provided that its overhead remains sufficiently small. This situation frequently arises in implementations of interpreted programming languages such as MATLAB: although the user-defined evaluation of the ODE right-hand side is performed in MATLAB, the linear-algebra operations responsible for most of the overhead are executed by highly optimized underlying routines~\cite{Rackauckas2017,Hairer2003}.

However, in other programming environments such as the Julia language, the 
work--precision diagrams based on function evaluations may not reflect the true 
computational efficiency.To quantify this, we introduce the CPU time 
per function evaluation ratio (CPU-FCN ratio for short), defined as the ratio of the total 
CPU time required for a given integration to the total number of function evaluations 
performed during that integration. For explicit splitting methods, this ratio is essentially 
constant: it depends neither on the step size nor on the particular scheme (i.e., the 
number of flows evaluated per step), since the computational cost is entirely dominated by 
the function evaluations themselves. For IRKGL16, by contrast, each step involves three 
additional sources of overhead beyond the function evaluations: the initial extrapolation 
to provide a starting guess for the fixed-point iteration, the overhead of the fixed-point 
iterations themselves, and the updating of state variables at the end of the step. 
Table~\ref{tab:overhead} reports, for each of the second-order ODE problems, the CPU-FCN 
ratio of IRKGL16 normalized by that of the explicit splitting methods, for both the fully 
sequential implementation IRKGL16-SEQ (\verb|simd=false, fseq=true|) and the fully 
vectorized implementation IRKGL16-SIMD (\verb|simd=true, fseq=false|). For IRKGL16-SEQ, 
the normalized ratio exceeds two in both test problems, indicating that the cumulative 
overhead surpasses the cost of the function evaluations themselves. For IRKGL16-SIMD, by 
contrast, the normalized ratio is well below one. This is because the SIMD 
vectorization not only reduces the overhead of the implicit solver, but also enables the 
8 stage function evaluations within each fixed-point iteration to be executed in parallel 
across SIMD registers, effectively amortizing their cost. As a consequence, each function 
evaluation in IRKGL16-SIMD is significantly cheaper in terms of wall-clock time than in 
the explicit methods, where function evaluations are performed sequentially. This implies that, 
relative to the work--precision diagrams presented above, the true efficiency curve of 
IRKGL16-SIMD is shifted to the left, making it even more competitive against the explicit 
symplectic integrators than those diagrams suggest.

\begin{table}[t]
\centering
\caption{CPU-FCN ratio of IRKGL16 normalized by that of the explicit splitting methods,
for the two second-order ODE test problems. ``6-body Solar System (vec)'' denotes the  SIMD-vectorized implementation of the flows }
\label{tab:overhead}
\begin{tabular}{lcc}
\hline
 & IRKGL16-SEQ & IRKGL16-SIMD \\
\hline
Hénon--Heiles                          & 2.5 & 0.6 \\
6-body Solar System                    & 2.6 & 0.6 \\
6-body Solar System (vec)        & 3.5 & 0.8 \\
\hline
\end{tabular}
\end{table}

It should be noted that the CPU-FCN ratio of IRKGL16 is not constant across step sizes. 
While the overhead of the fixed-point iterations scales proportionally with the number of 
function evaluations per step, the overhead of the initial extrapolation and the state 
update does not. Consequently, for smaller step sizes---which require fewer fixed-point 
iterations per step---these fixed costs carry more relative weight, leading to a slightly 
larger CPU-FCN ratio. The values reported in Table~\ref{tab:overhead} are computed at a 
reference step size for which the relative local energy error is approximately $10^{-14}$.

It is worth noting that explicit splitting methods can also benefit from SIMD 
vectorization, by explicitly exploiting SIMD registers in the implementation of the 
individual flows. To illustrate this, we have implemented a SIMD-optimized version of 
the flows for the 6-body outer Solar System problem, which results in a reduced CPU-FCN 
ratio for the explicit methods. We remark that SIMD-optimizing the flows requires 
additional implementation effort from the user, whereas the SIMD vectorization in IRKGL16 
is handled internally with a generic implementation of the user-defined function $g$.

To assess the impact of this 
optimization on the relative computational cost of IRKGL16, Table~\ref{tab:overhead} also reports 
(last row) the CPU-FCN ratio of IRKGL16 normalized by that of the SIMD-optimized explicit 
implementation. As expected, the normalized ratios increase with respect to the 
non-vectorized case, reflecting the fact that the explicit methods now perform each 
function evaluation more cheaply. Nevertheless, IRKGL16-SIMD retains a normalized ratio 
well below one, confirming that its fully vectorized implementation remains competitive 
even when compared against SIMD-optimized explicit integrators.

\begin{figure}[t]
		\centering
		\includegraphics[width=0.8 \textwidth]{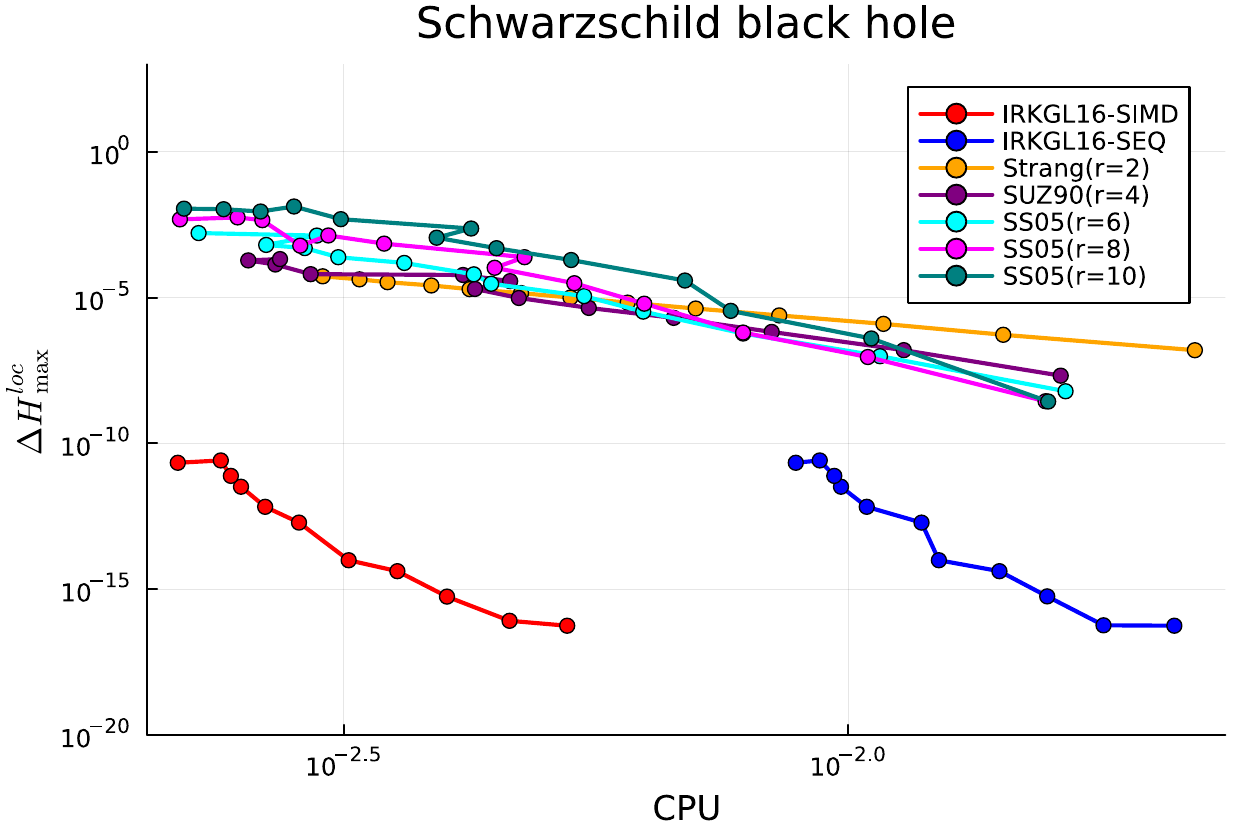} 
			\caption{Work-precision diagrams (maximum local Hamiltonian error versus CPU time) for the Schwarzschild black hole problem using symmetric composition schemes of orders $r = 2, 4, 6, 8$, and $10$, and two variants of IRKGL16: (i) fully vectorized (IRKGL16-SIMD), (ii) fully sequential (IRKGL16-SEQ)}
	\label{fig1}
\end{figure}

Figure~\ref{fig1} presents work--precision diagrams for the Schwarzschild black hole 
problem, displaying the maximum local error in the Hamiltonian \eqref{eq:DHlocmax} 
plotted against CPU time for explicit symplectic methods of orders $2,4,6,8,10$ and for 
two variants of IRKGL16: (i) the fully vectorized implementation IRKGL16-SIMD 
(\verb|simd=true, fseq=false|), and (ii) the fully sequential implementation IRKGL16-SEQ 
(\verb|simd=false, fseq=true|). In this example, we plot directly against CPU time rather 
than function evaluations, since the evaluation of the $t$-flows required to compute the 
basic map $\phi_h$ in \eqref{eq:genStrang} is considerably more expensive than the 
evaluation of the right-hand side of the ODE system, making the number of function 
evaluations a poor proxy for computational cost.

We observe that IRKGL16 outperforms the explicit symplectic integrators across the entire 
accuracy range shown, with IRKGL16-SEQ already surpassing them in the high-accuracy 
regime. A key factor is that the $t$-flow evaluations in \eqref{eq:genStrang} are 
considerably more expensive than the right-hand side evaluations in IRKGL16. As a result, 
for any given CPU-time budget, IRKGL16-SIMD delivers substantially higher accuracy than 
the explicit methods, regardless of the accuracy level they achieve.

\subsubsection{Evolution of errors for the high-precision regime}

\begin{figure}[t]
	\centering	
	\begin{subfigure}{0.48\textwidth}
		\centering
		\includegraphics[width=\textwidth]{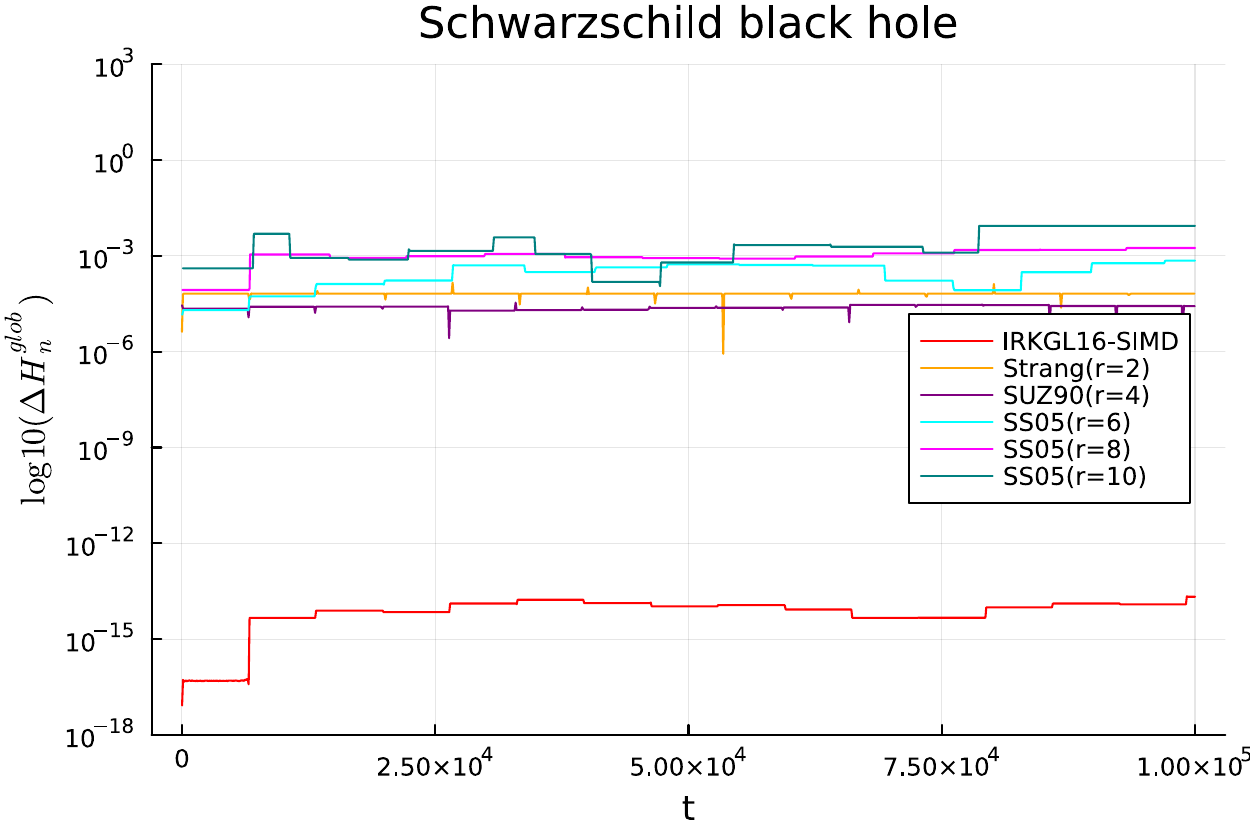}
		\caption{}
		\label{fig:Fig2a}
	\end{subfigure}%
	\hfill
	\begin{subfigure}{0.48\textwidth}
		\centering
		\includegraphics[width=\textwidth]{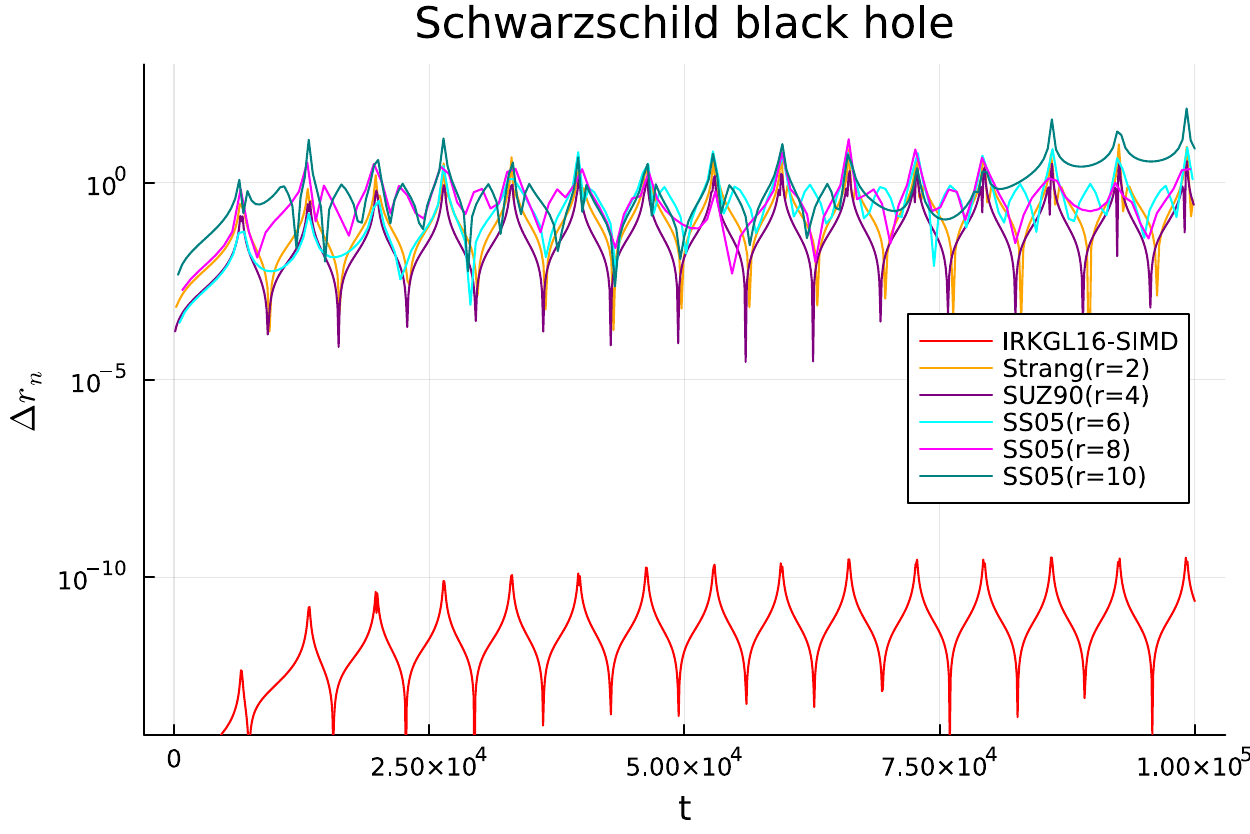}
		\caption{}
		\label{fig:Fig2b}
	\end{subfigure}%
	\caption{Left: global error in the Hamiltonian. Right: relative error in the radial coordinate $r$ for the Schwarzschild black hole problem. Results are shown for IRKGL16-SIMD and symmetric composition methods of orders $r = 2, 4, 6, 8, 10$}
	\label{fig2}
\end{figure}

We now examine the time evolution of both the global Hamiltonian error and the position
error for each of the three test problems in the high-precision regime. For each problem,
we select a reference step size $h_{\mathrm{ref}}$ for IRKGL16-SIMD such that the maximum
global Hamiltonian error
\begin{align}
	\label{eq:DHglobmax}
	\Delta H^{\text{glob}}_{\max} = \max_n \left( \Delta H^{\text{glob}}_n \right), \quad \Delta H^{\text{glob}}_n = \left|\frac{H(y_n) - H(y_0)}{H(y_0)}\right|,
	\end{align}
 is approximately $10^{-14}$. Each explicit
method is then run with a step size chosen so that its total CPU time matches that of
IRKGL16-SIMD at $h_{\mathrm{ref}}$. This ensures a fair comparison on equal computational
budget grounds.

For the Schwarzschild black hole problem, the results are shown in Figure~\ref{fig2}. The
left panel displays the global Hamiltonian error and the right panel the relative error in
the radial coordinate $r$, defined as $\Delta r_n = |r_n - \tilde{r}_n|/|\tilde{r}_n|$,
where $\tilde{r}_n$ denotes the reference solution computed with IRKGL16-SIMD at a
smaller step size ($h=h_{\mathrm{ref}}/4$). We compare IRKGL16-SIMD against symmetric composition methods of orders
$r = 2, 4, 6, 8, 10$. At the CPU-time-matched step sizes, the most competitive explicit
methods are those of orders 2 and 4, yet their global Hamiltonian and position errors
remain approximately nine orders of magnitude larger than those of IRKGL16-SIMD throughout
the integration.

\begin{figure}[t]
	\centering	
	\begin{subfigure}{0.48\textwidth}
		\centering
		\includegraphics[width=\textwidth]{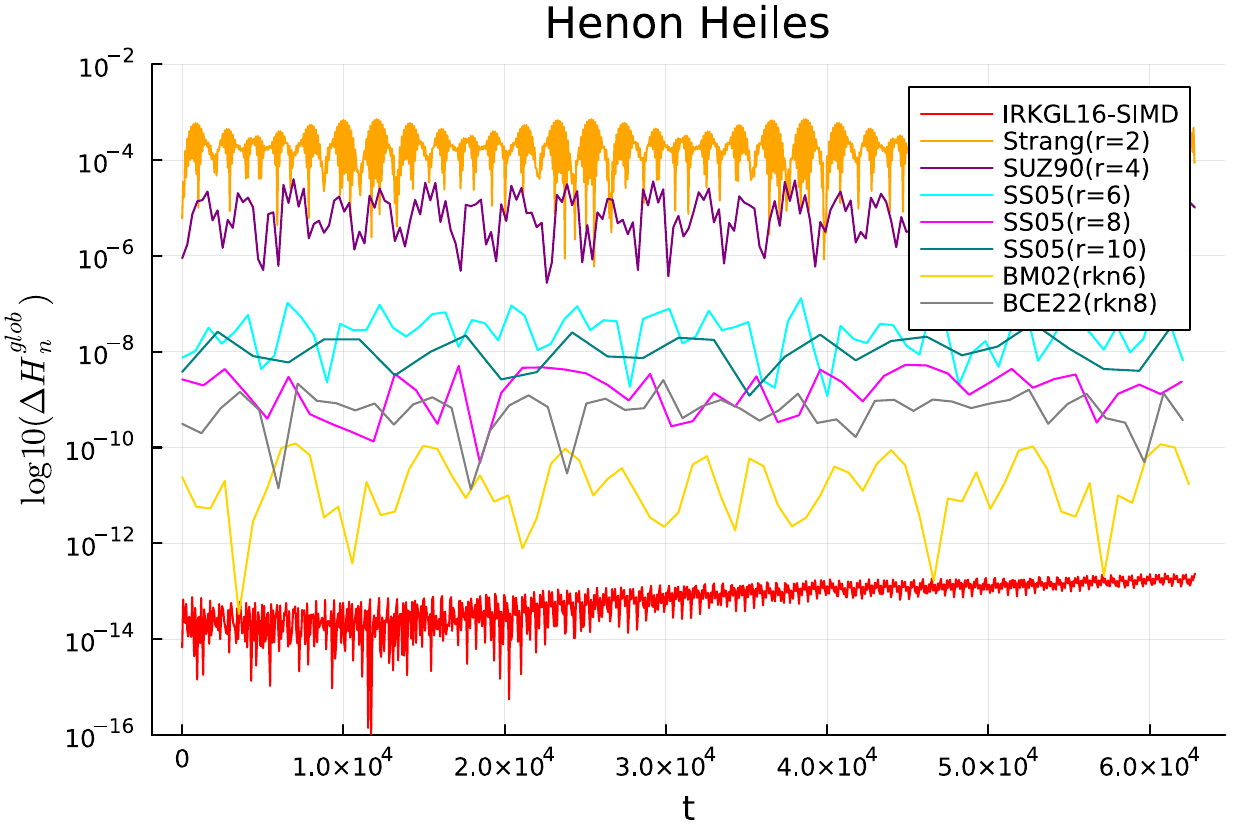}
		\caption{}
		\label{fig:Fig4a}
	\end{subfigure}%
	\hfill
	\begin{subfigure}{0.48\textwidth}
		\centering
		\includegraphics[width=\textwidth]{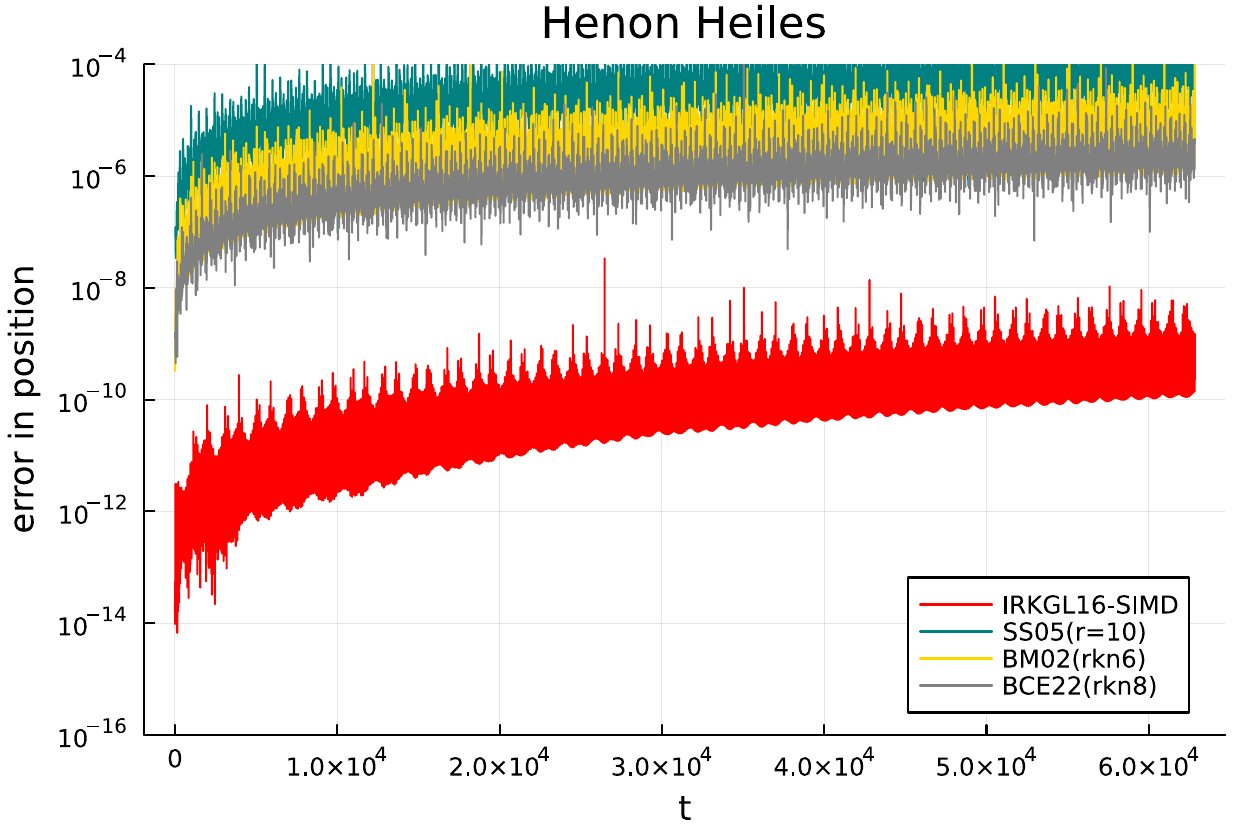}
		\caption{}
		\label{fig:Fig4b}
	\end{subfigure}%
	\caption{Hénon-Heiles problem: time evolution of (i) the global Hamiltonian error (left), 
	(ii) error in position (right)  using IRKGL16, symmetric compositions, and RKN methods}
	\label{fig4}
\end{figure}

Figure~\ref{fig4} presents the analogous comparison for the H\'{e}non-Heiles problem,
showing the global Hamiltonian error (left) and the absolute error in position $q$ (right),
with the reference solution again provided by IRKGL16-SIMD at small step size. Here the
most competitive explicit methods at equal CPU cost are the RKN-type integrators BM02
(order 6) and BCE22 (order 8), yet IRKGL16-SIMD maintains a substantially lower error in
both the Hamiltonian and position throughout the integration interval.

\begin{figure}[h!]
	\centering	
		\includegraphics[width=0.8\textwidth]{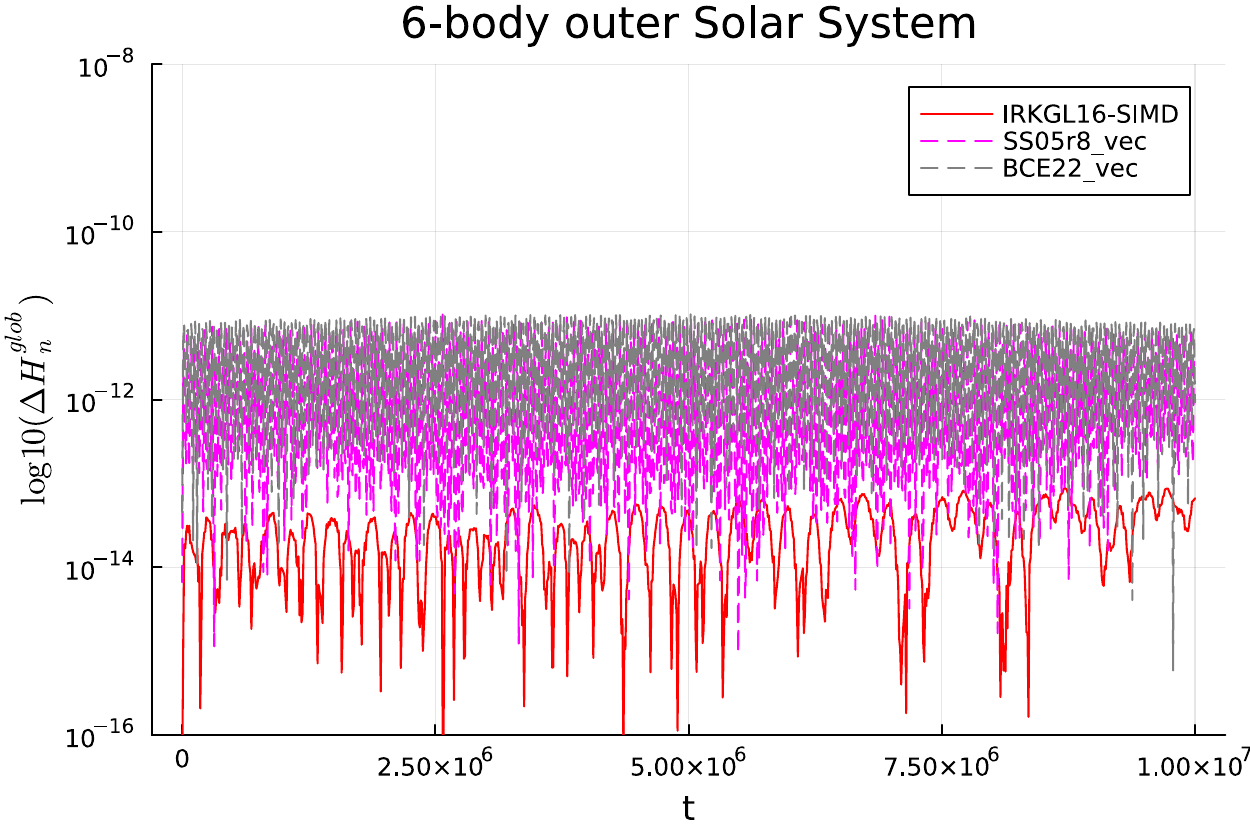}
%	\caption{Time evolution of the global Hamiltonian error for the 6-body outer Solar System for IRKGL16, and the SIMD-optimized version of the 8th order integrators SS05r8 and  BCE22}
	\caption{Time evolution of the global Hamiltonian error for the outer Solar System 6-body problem, computed with IRKGL16 and the 8th-order integrators SS05r8 and BCE22 using SIMD-vectorization of the flows.}
	\label{fig5}
\end{figure}

\begin{figure}[h!]
	\centering	
	\begin{subfigure}{0.48\textwidth}
		\centering
		\includegraphics[width=\textwidth]{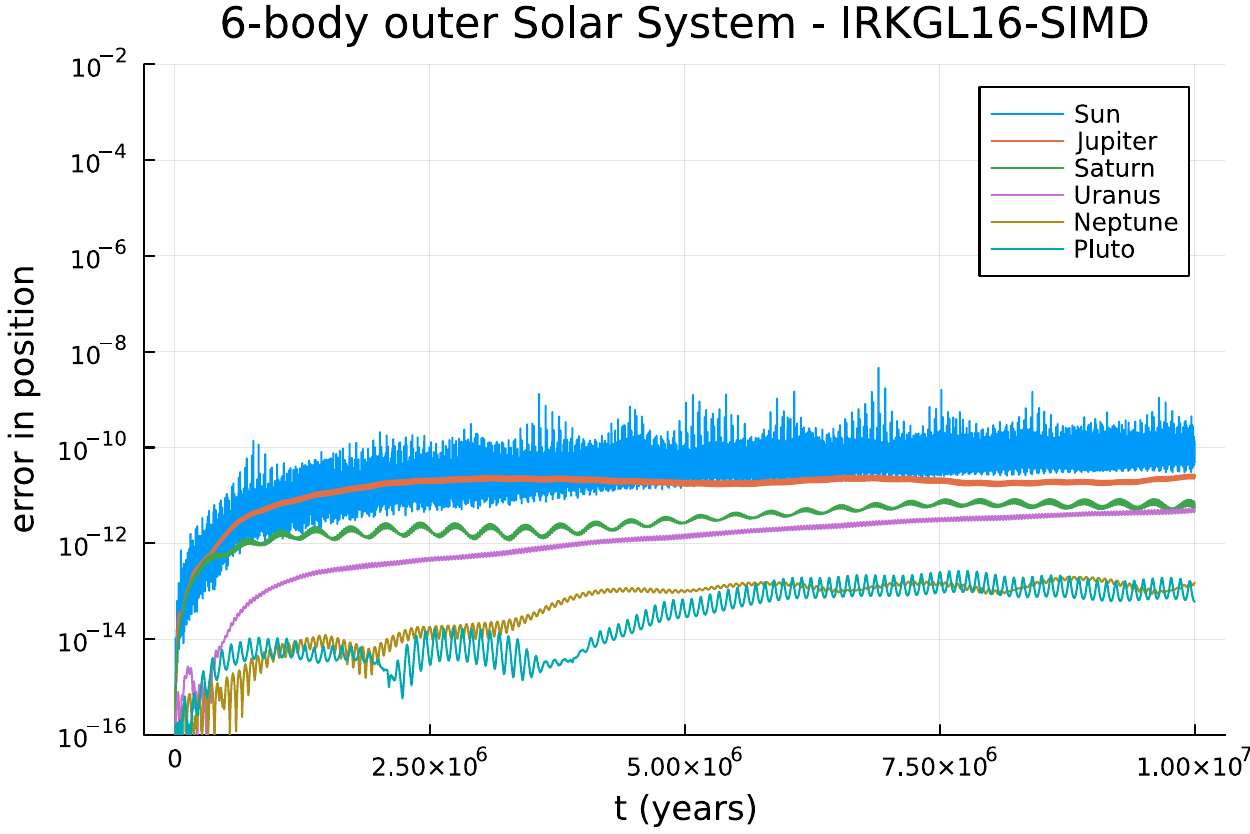}
		\caption{}
		\label{fig:Fig6a}
	\end{subfigure}%
	\hfill
	\begin{subfigure}{0.48\textwidth}
		\centering
		\includegraphics[width=\textwidth]{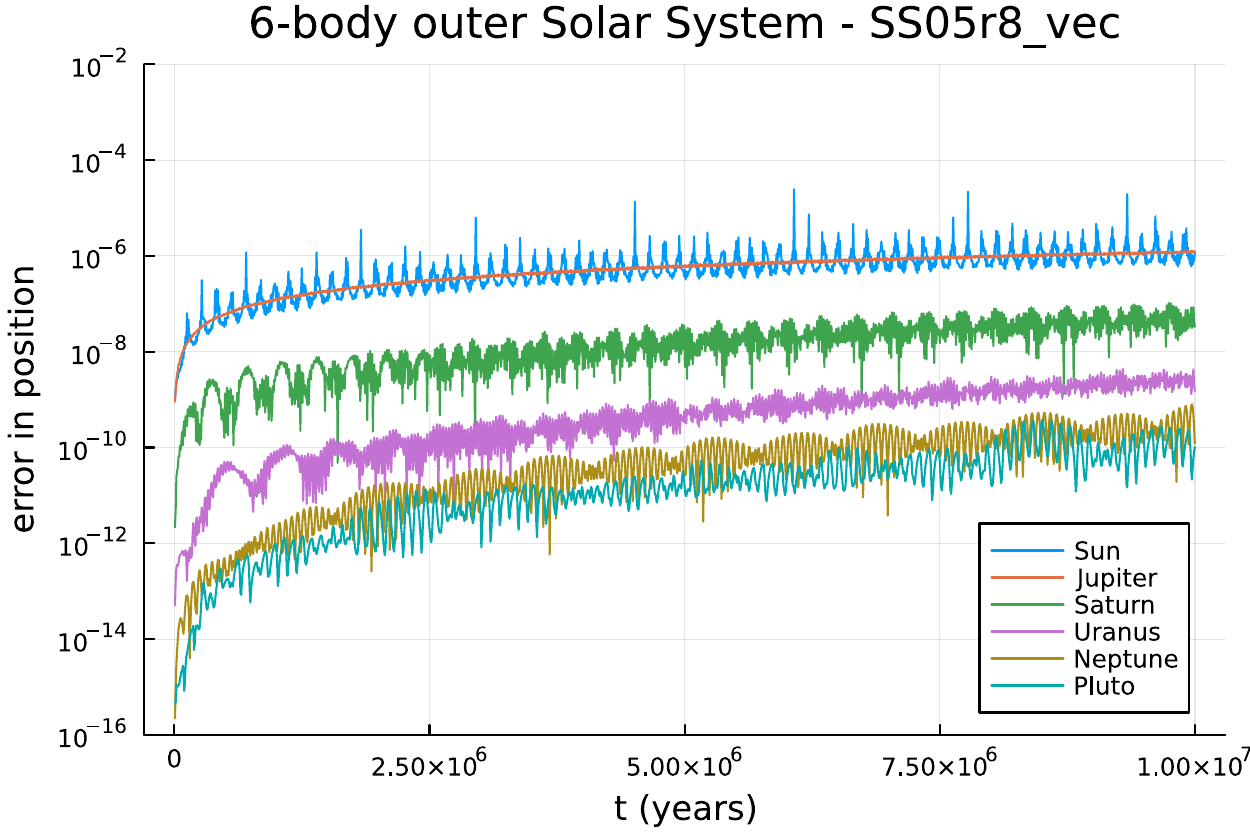}
		\caption{}
		\label{fig:Fig6b}
	\end{subfigure}%
%	\caption{Evolution of position error for each planet in the 6-body outer Solar System using IRKGL16 (left) and 
%	the SIMD-optimized version of the 8th order integrator SS05r8}

	\caption{Evolution of the position error for each planet in the outer Solar System 6-body problem: IRKGL16 (left) and the 8th-order integrator SS05r8  using SIMD-vectorization of the flows (right).}
	\label{fig6}
\end{figure}

For the 6-body outer Solar System, the comparison is carried out against SIMD-optimized
implementations of the two most competitive 8th-order explicit integrators, SS05r8 and
BCE22, with step sizes matched to the CPU time of IRKGL16-SIMD using the SIMD-optimized
evaluation of the acceleration function $g$. This represents a particularly challenging
benchmark for IRKGL16, as the explicit methods benefit from problem-specific SIMD
acceleration. Figure~\ref{fig5} shows the time evolution of the global Hamiltonian error:
even under these conditions, IRKGL16-SIMD achieves lower Hamiltonian errors than both
explicit competitors throughout the integration. The evolution of position errors for each
planet is shown in Figure~\ref{fig6}, comparing IRKGL16-SIMD (left) against SS05r8
(right), which exhibits lower position errors than BCE22; again, IRKGL16-SIMD yields lower position errors across all planets over the
entire integration interval.

Across all three test problems considered---the Schwarzschild black hole, the
H\'{e}non-Heiles system, and the 6-body outer Solar System---IRKGL16-SIMD consistently
achieves lower global Hamiltonian and position errors than the best-performing explicit
symplectic integrators at equal CPU time in the high-precision regime.

\subsection{Discussion on alternative implementations}

The main source of overhead in implicit methods implemented with fixed-point iterations
arises from the repeated updating of the state variables stored in the $D$-dimensional
arrays $U_i$ (and $L_i$) at each iteration. At least two distinct approaches can be adopted to
exploit SIMD registers in order to reduce this overhead.

The first approach performs the updates of each state variable in the $D$-dimensional
arrays $U_i$ independently, and vectorizes them component-wise, provided that the state
components of each array $U_i$ are stored consecutively in memory. This is the more
standard approach, and can in principle be handled automatically by the compiler, provided
that the code is written appropriately. Specifically, the inner loops over $j = 1, \dots,
s$ required for the stage summations \eqref{eq:Yi2} or \eqref{eq:Qi2} must be unrolled, so that the outer loop over $k = 1,
\dots, D$ consists of straightforward arithmetic that the compiler can SIMD-vectorize. We
have verified that, in our sequential implementation of IRKGL16, unrolling these inner
loops and disabling bounds checking via the \verb|@inbounds| macro causes the Julia
compiler to apply component-wise vectorization automatically, substantially reducing the
overhead. A further, marginal reduction can be achieved by applying, instead of
\verb|@inbounds|, the \verb|@turbo| macro from the \texttt{LoopVectorization.jl} package.

The second approach instead arranges, for each state variable, the corresponding
components $U_1, \dots, U_s$ consecutively in memory, and vectorizes their simultaneous
update using SIMD $s$-vectors. This is the approach adopted in our implementation of
IRKGL16, explained in detail in Subsection~2.5 and Section~3. We have verified that when
the right-hand side function evaluations are performed sequentially, both vectorization
approaches---the stage-wise vectorization of IRKGL16 with \texttt{simd=true, fseq=false}
on one hand, and the component-wise vectorization via unrolled inner loops with the
\verb|@turbo| macro on the other---yield essentially the same reduction in overhead. The
key advantage of stage-wise vectorization is that it enables the $s$ evaluations of
the right-hand side function to be performed in parallel within SIMD registers, whenever
it is implemented as a generic function compatible with the \texttt{SIMD.jl} package. This
SIMD-parallel evaluation is handled internally by IRKGL16 and remains completely
transparent to the user, who is only required to provide a standard generic implementation
of the right-hand side function. In summary, we advocate for stage-wise vectorization
---and have adopted it in IRKGL16 in particular---as a convenient and user-friendly
implementation strategy for IRKGL schemes.

To further assess the practical relevance of component-wise vectorization, we computed the
CPU-FCN ratios of the \texttt{Fortran~77} implementations of both the 12th-order IRKGL
scheme and the explicit composition schemes provided in the \texttt{gnicodes} folder of
the software repository accompanying~\cite{Hairer2006}, applied to the 6-body outer Solar
System problem and compiled with the \verb|-O2| and \verb|-O3| optimization flags, the latter requesting vectorization whenever possible. We do
not consider the \verb|-Ofast| flag, since it enables aggressive optimization that silently eliminates compensated summation algorithms by treating
their error-compensation steps as algebraically redundant, 
thereby compromising the long-term energy error propagation control that symplectic 
integrators rely on.
For the explicit composition schemes, the CPU-FCN ratio normalized by that of our
SIMD-optimized Julia implementation of the flows is approximately 1.0 under
\verb|-O3| and 1.2 under \verb|-O2|, suggesting that the compiler successfully applies
SIMD acceleration under \verb|-O3|. This also shows that similar efficiency can be
achieved for explicit symplectic schemes in both Fortran and Julia, provided that the
implementation is set up so that the compiler can exploit SIMD acceleration. For IRKGL12,
by contrast, the normalized CPU-FCN ratio is approximately 1.4 under both \verb|-O2| and
\verb|-O3|, indicating that the compiler fails to apply vectorization even under \verb|-O3|.
 The reason is that the number of stages $s$ is treated as a
generic parameter in the \texttt{gnicodes} implementation, which prevents the inner
summation loops from being unrolled. We have verified that setting $s = 6$ statically, so
that the number of terms in the summation is known at compile time, allows the compiler to
unroll the loop and successfully apply component-wise vectorization under \verb|-O3|,
reducing the normalized CPU-FCN ratio to approximately 1.05 (compared to 1.4 under
\verb|-O2|). This confirms that automatic component-wise vectorization is also achievable
in Fortran implementations of IRKGL schemes, provided the number of stages is fixed at
compile time. However, even in this favorable case, the normalized CPU-FCN ratio remains
above one, since no parallelism is exploited across the $s$ independent evaluations of the
right-hand side function. By contrast, IRKGL16-SIMD achieves a normalized CPU-FCN ratio
of approximately 0.8 in the 6-body Solar System problem (normalized against explicit
schemes with SIMD-optimized right-hand-side function), by seamlessly and transparently
exploiting stage-wise SIMD parallelism across all function evaluations.

For ODE systems whose right-hand side is not compatible with \texttt{SIMD.jl}, an
alternative strategy is to evaluate the $s$ stage functions in parallel using
multithreading. In Julia, this can be conveniently implemented using the
\texttt{Threads.@threads} macro. However, the additional synchronization overhead incurred
by thread-level parallelism is only compensated when the right-hand side function is
sufficiently expensive to evaluate \cite{Hairer2008}. We have verified that, for the three test problems
considered in this work, a straightforward application of \texttt{Threads.@threads} to
parallelize the stage function evaluations results in increased total CPU time, indicating
that the right-hand side functions in these examples are not costly enough to amortize the
threading overhead.

\section{Conclusions}\label{sec13}

High-order explicit symplectic integrators are effective when the Hamiltonian can be
decomposed into solvable components, but IRKGL-based integrators offer a powerful
alternative: they are symplectic, time-reversible, and high-order accurate, can be
efficiently implemented using fixed-point iterations in non-stiff regimes, and are broadly
applicable to general canonical Hamiltonian systems~\cite{Sanz1994,Hairer2006} in
first-order ODE form, regardless of whether the Hamiltonian can be split into solvable
parts.

In Section~2, we presented a general stage-wise vectorized formulation of IRKGL schemes
for both first-order and second-order Hamiltonian systems. Beyond reducing computational
overhead, this formulation leaves the algorithm naturally ready for stage-wise SIMD
vectorization, enabling parallel evaluation of the right-hand side function across all
stages whenever the function is compatible with \texttt{SIMD.jl}. For second-order
systems, we introduced a reformulation of the corresponding Runge--Kutta--Nyström method
that ensures exact symplecticity at the level of double-precision floating-point
arithmetic, extending the approach of~\cite{antonana2017} originally developed for
first-order systems to this setting.

The numerical experiments of Section~5 confirm the effectiveness of the proposed approach.
For Hamiltonian problems admitting a decomposition into a sum of exactly solvable components, IRKGL16 outperforms optimized
high-order explicit symplectic integrators in the high-precision regime across all test
problems considered. When the right-hand
side function is compatible with \texttt{SIMD.jl}, this advantage is most pronounced; when
it is not, the stage-wise vectorization of all remaining computations still greatly reduces
overhead, and IRKGL16 may still surpass explicit methods at high-precision levels.

Overall, our results show that IRKGL integrators, optimized with the proposed
stage-wise SIMD vectorization technique, provide a powerful and general alternative to
explicit symplectic methods, achieving high accuracy efficiently for non-stiff Hamiltonian
systems.

\backmatter

\bmhead{Acknowledgments}

All the authors have received funding by the Spanish State Research Agency through project PID2022-136585NB-C22, MCIN/AEI/10.13039/501100011033 and the European Union. They are also partially supported by the Department of Education of the Basque Government through the Consolidated Research Group MATHMODE (ITI456-22).

%\section*{Nire bibliografia}


%\bibliography{bibliography}% common bib file

\bibliographystyle{spmpsci}      % mathematics and physical sciences

\end{document}